\documentclass[11pt, oneside]{amsart}
\usepackage[letterpaper,margin=1in,heightrounded]{geometry}
\usepackage[utf8]{inputenc}
\usepackage{amsfonts, amsmath, amssymb, amsthm, bm}
\usepackage{aliascnt, hyperref, cleveref}
\hypersetup{hidelinks,
  pdftitle={Height function covariance and point-to-line coalescence in the directed landscape},
  pdfauthor={Yu Gu, Jeremy Quastel, Tianhao Xian},
  pdfkeywords={Directed landscape, KPZ fixed point, geodesic coalescence, integration by parts, Brownian last-passage percolation, O'Connell-Yor polymer, Airy 1 process}}
\usepackage{xcolor}
\usepackage{graphicx}
\usepackage{longtable}
\usepackage{enumerate}
\usepackage[numbers,sort&compress]{natbib}

\newcommand{\jbrac}[1]{\left\langle {#1}\right\rangle}
\newcommand{\pare}[1]{\left({#1}\right)}
\newcommand{\sqbrac}[1]{\left[ {#1}\right]}
\newcommand{\brac}[1]{\left\{{#1}\right\}}

\newcommand{\R}{\mathbb{R}}
\newcommand{\E}{\mathbb{E}}
\newcommand{\Z}{\mathbb{Z}}
\renewcommand{\P}{\mathbb{P}}
\newcommand{\bE}{\mathbf{E}}
\newcommand{\bP}{\mathbf{P}}
\newcommand{\bs}{\mathbf{s}}
\newcommand{\bt}{\mathbf{t}}
\newcommand{\cD}{\mathcal{D}}
\newcommand{\cL}{\mathcal{L}}
\newcommand{\cS}{\mathcal{S}}
\newcommand{\fh}{\mathfrak{h}}
\newcommand{\fu}{\mathfrak{u}}

\newcommand{\ind}[1]{\mathbf{1}_{#1}}

\newcommand{\Epoly}{\E^{\mathrm{poly},\beta}}
\newcommand{\Ppoly}{\P^{\mathrm{poly},\beta}}

\newcommand{\sgn}{\operatorname{sgn}}

\DeclareMathOperator{\cov}{Cov}
\DeclareMathOperator{\Var}{Var}

\reversemarginpar
\newtheorem{theorem}{Theorem}[section]
\newtheorem{proposition}[theorem]{Proposition}
\newtheorem{lemma}[theorem]{Lemma}

\theoremstyle{definition}

\theoremstyle{remark}

\title[Height covariance and point-to-line coalescence]{Height function covariance and point-to-line coalescence in the directed landscape}
\author{Yu Gu}
\address[Yu Gu]{Department of Mathematics, University of Maryland, College Park, Maryland, USA}
\email{ygu7@umd.edu}

\author{Jeremy Quastel}
\address[Jeremy Quastel]{Department of Mathematics, University of Toronto, Toronto, Ontario, Canada}
\email{quastel@math.toronto.edu}

\author{Tianhao Xian}
\address[Tianhao Xian]{Department of Mathematics, University of Toronto, Toronto, Ontario, Canada}
\email{tianhao.xian@utoronto.ca}

\date{June 2026}

\begin{document}

\begin{abstract}
We prove an integration-by-parts formula for the KPZ fixed point with a deterministic $C_b^1$ terminal reward. For flat data, let $C(x)$ denote the spatial height covariance at time one. For almost every $x>0$, we show that $-C'(x)$ is the probability that the point-to-line geodesics started at $0$ and $x$ coalesce by time one. The proof uses Gaussian integration by parts for the O'Connell--Yor polymer, followed by the zero-temperature and directed-landscape limits.
\end{abstract}

\maketitle
\tableofcontents

\section{Introduction}

The directed landscape is the universal scaling limit of last-passage percolation in the KPZ universality class \cite{dauvergne2022directed}. Besides describing fluctuations of passage times, it carries a random geometry of optimizing paths, or geodesics. A fundamental feature of this geometry is \emph{coalescence}: geodesics starting at different spatial points may meet and thereafter follow the same path. The purpose of this paper is to relate this geometric event to the spatial covariance of the KPZ fixed point, and to establish an integration-by-parts identity which explains the relation.

Let $\cL(x,s;y,t)$ be the directed landscape. For a terminal reward $\fh_0$, consider the point-to-line variational problem
\begin{equation}\label{eq:introheight}
\fh(x)=\sup_{y\in\R}\{\cL(x,0;y,1)+\fh_0(y)\}.
\end{equation}
This is the $t=1$ KPZ fixed-point height function. Write $\gamma_t(x)$ for the rightmost optimizing path from $(x,0)$ to its maximizing terminal point at time $1$. For two starting points $x<y$, let $\tau_{x,y}$ be the first time their paths meet, with $\tau_{x,y}=\infty$ if they do not meet by time $1$. Once two optimizing paths meet, they coalesce. Thus $\{\tau_{x,y}\le1\}$ is the event that their optimizing paths have joined during the unit time interval.

For flat reward $\fh_0\equiv0$, the process $\fh$ is spatially stationary. In our normalization, it is related to the standard Airy$_1$ process by $\fh(x)\stackrel{\rm d}{=}2^{1/3}\mathcal A_1(2^{-2/3}x)$, as processes. Set
\[
C(x)=\cov(\fh(0),\fh(x)),\qquad
p(x)=\bP(\tau_{0,x}\le1),\qquad x\ge0.
\]
Our first main result is the exact identity
\begin{equation}\label{eq:introcov}
C(a)-C(b)=\int_a^b p(x)\,dx,\qquad 0\le a<b.
\end{equation}
In particular, $C$ is absolutely continuous and $-C'(x)=p(x)$ for almost every $x>0$. This identifies a two-geodesic coalescence probability with the spatial derivative of a two-point height covariance.

The geometric meaning of the identity can be expressed through the endpoint map. Under flat reward, let $Y(x)=\gamma_1(x)$ be the rightmost maximizing endpoint. The map $x\mapsto Y(x)$ is nondecreasing and has intervals of constancy. Its jump points form the point-to-line \emph{lower point field}, also called the shock-point field; see, for example, \cite{rahman2025infinite,dauvergne2023last,dauvergne2025geodesic,bakhtin2018global}. They separate starting points whose optimizing geodesics reach different terminal endpoints. For $x<y$, the event that the two geodesics coalesce by time $1$ is, up to conventions at the endpoints, the event that there is no separating jump in $(x,y)$. Thus $p(x)$ is a \emph{void probability} for this random point field. Formula~\eqref{eq:introcov} identifies this void probability with the negative derivative of the Airy$_1$ covariance.

Let $\mathcal P$ denote the lower point field, and let
\[
D=\inf\{x>0:x\in\mathcal P\}
\]
be the distance from the origin to its next point on the right. The void-probability interpretation gives $\bP(D>x)=p(x)$, and hence, using the decay of the Airy$_1$ covariance at infinity,
\begin{equation}\label{eq:intromeandistance}
\bE D=\int_0^\infty p(x)\,dx=C(0)=\Var(\fh(0)).
\end{equation}
Thus the mean distance to the next shock point equals the one-point height variance. Moreover, the cubic large-distance decay discussed below implies that $D$ has moments of every order. More precisely, for every $k>0$,
\[
\bE D^k=k\int_0^\infty x^{k-1}p(x)\,dx.
\]
For $k>1$, integration by parts in \eqref{eq:introcov} gives
\begin{equation}\label{eq:introDmoms}
\bE D^k=k(k-1)\int_0^\infty x^{k-2}C(x)\,dx,
\qquad\text{in particular}\qquad
\bE D^2=2\int_0^\infty C(x)\,dx.
\end{equation}

Another equivalent form of \eqref{eq:introcov}, worth recording, is the height-increment identity
\begin{equation}\label{eq:introvariance}
\bE[(\fh(L)-\fh(0))^2]
=2\int_0^L p(x)\,dx=2\bE(D\wedge L).
\end{equation}
It can also be obtained directly from the integration-by-parts formula below by formally taking $g=\mathbf1_{[0,L]}$ and $F(\fh)=\langle g,\partial_x\fh\rangle$, and then approximating these test functions. In particular, the left-hand side tends to $2\bE D=2\Var(\fh(0))$ as $L\to\infty$.

\paragraph{The flat intensity.}
The Airy$_1$ boundary-layer estimate proved in Appendix~\ref{app:flat-curvature}
gives the short-distance expansion
\begin{equation}\label{eq:flatcurvatureintro}
\bE[(\fh(x)-\fh(0))^2]
=2x+\frac43\bE\fh(0)\,x^2+o(x^2),\qquad x\downarrow0.
\end{equation}
The covariance--coalescence identity and monotonicity of the void probability then give
\begin{equation}\label{eq:flat-intensity-candidate}
\rho:=\lim_{x\downarrow0}\frac{1-p(x)}x
=-\frac43\bE\fh(0)
=-\frac{2^{4/3}}3\bE\chi_{\rm GOE} \approx 1.0134.
\end{equation}
Here $\chi_{\rm GOE}$ has the standard GOE Tracy--Widom law. The proof,
including the third-order uniform Airy$_1$ Fredholm-minor remainder, is given
in Appendix~\ref{app:flat-curvature}.

For $\fh_0(y)=\sigma B(y)$, where $B$ is two-sided Brownian
motion with variance $1$ per unit distance, independent of the directed landscape, we also have
\[
V_\sigma(t,L)=\bE\big[(\fh_\sigma(t,L)-\fh_\sigma(t,0))^2\big],\qquad
p_\sigma(t,L)=\bP_\sigma(\tau_{0,L}\le t).
\]
In the variance-$2$ Brownian normalization of the KPZ fixed point,
the corresponding Gaussian integration-by-parts relation takes the form
\begin{equation}\label{eq:intro-sigma-brownian}
\frac12\partial_L V_\sigma(t,L)
=\frac{\sigma^2}{2}
+\left(1-\frac{\sigma^2}{2}\right)p_\sigma(t,L),
\quad L>0,
\end{equation}
or, equivalently,
\[
V_\sigma(t,L)=\sigma^2L+(2-\sigma^2)\int_0^L p_\sigma(t,x)\,dx.
\]
For $\sigma=0$ this is \eqref{eq:introvariance}; at $\sigma^2=2$ it
reduces to the Brownian invariant increment variance $2L$.
The coefficient $1-\sigma^2/2$ arises from the additional Gaussian
integration by parts in the random initial profile, which cancels part
of the environmental coalescence contribution; compare
\cite{gu2024integration,pimentel2022integration}.
The proof is the same as for flat except that we use Prop. 2.3 of \cite{gu2024integration} to integrate the term in \eqref{eq:introibp} by parts.  For general $\sigma$, the
natural formulation is in terms of the \emph{increment variance} $V_\sigma$:
one should not replace its derivative by that of the height covariance,
since the unanchored heights are not spatially stationary.

The relationship between correlations and coalescence has a long history in interacting particle systems; see, for example, the classical coalescing-walk construction of Arratia \cite{arratia1979coalescing}. A particularly transparent example is the duality between the voter model and coalescing random walks \cite{holley1975ergodic}. If the voter model starts from independent Bernoulli$(q)$ spins, then the covariance of two spins at time $t$ is $q(1-q)$ times the probability that their backward ancestral walks have coalesced by time $t$. Correlations can therefore be understood through shared ancestry. Our identity is reminiscent of this principle, but it is different in two essential respects: a \emph{spatial derivative} of the covariance appears, and the paths are optimizing geodesics determined by a common random environment, not independent Markov paths that coalesce upon meeting.

Related connections have appeared in last-passage percolation, where coalescing geodesics are closely tied to exit points, competition interfaces, Busemann functions, and stationary measures. Relevant precedents include the stationary Hammersley-process work of Cator and Groeneboom \cite{cator2005hammersley}, Pimentel's duality between coalescence and exit points \cite{pimentel2016duality}, and the work of Georgiou, Rassoul-Agha, and Sepp\"al\"ainen on geodesics and competition interfaces \cite{georgiou2017geodesics}. Quantitative bounds for coalescence in related integrable polymer models have also been obtained by Rassoul-Agha, Sepp\"al\"ainen, and Shen \cite[Corollary~2.6]{rassoulagha2024coalescence}. Those settings often concern stationary or semi-infinite geodesics. Here we obtain an exact relation directly for the finite-time, flat-data KPZ fixed point: rather than bounding the covariance by a coalescence event or representing it through a different exit-point object, its derivative is the coalescence probability itself.

The mechanism behind \eqref{eq:introcov} is a more general integration-by-parts formula for the directed landscape. For $\fh_0\in C_c^1(\R)$, let $\fu=\partial_x\fh$ in the distributional sense. If $F$ is a bounded continuously differentiable cylindrical functional of the height, with variational derivative $D_yF$, then for $g\in C_c^1(\R)$ we prove
\begin{align}\label{eq:introibp}
\bE\bigl[F(\fh)\langle g,\fu\rangle\bigr]
&=\iint_{\R^2}g(x)\,\bE\bigl[D_yF(\fh)\,\sgn(y-x)\,\ind{\tau_{x,y}\le1}\bigr] \,dx\,dy\\
&\quad+\int_\R g(x)\,\bE\bigl[F(\fh)\fh_0'(\gamma_1(x))\bigr] \,dx.\notag
\end{align}
The precise class of test functionals is specified in Theorem~\ref{thm:KPZ fixed point integration by parts}. The first term involves the coalescence of two optimizing paths, while the second records the derivative of the terminal reward evaluated at the maximizing endpoint. For flat reward the latter disappears. Taking $F$ to approximate a height evaluation and $g$ to approximate an interval indicator gives \eqref{eq:introcov} after spatial stationarity is used.

Our proof begins with the O'Connell--Yor polymer in a Brownian environment. Gaussian integration by parts produces an overlap term involving two polymer paths. An exchange of the path tails at their first intersection cancels the contributions after that intersection, leaving precisely the probability of an intersection, with its orientation sign \cite[Theorem~1.1 and Corollary~1.3]{rassoulagha2026coalescence}. We then take the zero-temperature limit to Brownian last-passage percolation and the KPZ scaling limit to the directed landscape. The passage to the limit requires joint control of the height field, optimizing paths, and their first coalescence times. This route yields \eqref{eq:introibp} without appealing to an explicit Fredholm determinant formula for the KPZ fixed point.

\subsection*{Coalescence asymptotics}

The exact identity also makes it possible to transfer quantitative information about the Airy$_1$ process to the geometry of geodesics. Basu, Busani, and Ferrari \cite[Theorem~1.1]{basu2023exponent} proved sharp bounds on the cubic exponential scale of the Airy$_1$ covariance. In the normalization of the directed landscape used here, their estimates take the form, for sufficiently large $x$ and suitable positive constants,
\begin{equation}\label{eq:introairybound}
 e^{-x^3/3-cx\log x}\le C(x)\le e^{-x^3/3+c'x^2}.
\end{equation}
Since $p(x)$ is nonincreasing, the integrated covariance identity \eqref{eq:introcov} allows one to bound $p(x)$ between covariance differences on adjacent intervals. In particular,
\begin{equation}\label{eq:introtail}
 \lim_{x\to\infty}\frac1{x^3}\log\bP(\tau_{0,x}\le1)=-\frac13.
\end{equation}
Thus the distance $D$ to the next point of the lower point field has a cubic-exponential tail.

Rassoul-Agha, Sepp\"al\"ainen, and Shen \cite[Corollary~2.6]{rassoulagha2024coalescence} proved quantitative coalescence bounds for equilibrium log-gamma polymers. In the notation of the earlier version of this manuscript, their estimates were recorded as follows: for $\delta>c|x|^{-1/2}$,
\begin{equation}\label{eq:introRASS}
 e^{-c_1\delta^2}\le \bP(\tau_{0,x}<\delta |x|^{3/2})\le e^{-c_2\delta^2}.
\end{equation}
They also established the large-time estimates, for sufficiently large $t$,
\begin{equation}\label{eq:introRASSlong}
 \frac{c_1}{t^{2/3}}\le \bP(\tau_{0,1}>t)\le \frac{c_2(\log t)^{10}}{t^{2/3}}.
\end{equation}
Here the coalescence times refer to their equilibrium polymer geometry, rather than the finite-time point-to-line geodesics studied here.\footnote{The displayed short-time bounds are reproduced from the earlier draft and should be checked against the precise theorem and parameter convention of \cite[Corollary~2.6]{rassoulagha2024coalescence}: for small $\delta$, an exponent proportional to $-1/\delta^2$, rather than $-\delta^2$, is consistent with the small-time scaling discussed below.} To compare the scales, the directed landscape has the $1{:}2{:}3$ scaling relation
\begin{equation}\label{eq:introscaling}
 \tau_{0,ax}\ \stackrel{\rm d}{=}\ a^{3/2}\tau_{0,x},
\end{equation}
when the geodesics are defined on the whole positive time axis in a scaling-consistent way. Consequently, for $\delta>0$,
\begin{equation}\label{eq:introdelta}
 \bP(\tau_{0,x}<\delta x^{3/2})
 =\bP(\tau_{0,\delta^{-2/3}}<1).
\end{equation}
Our cubic tail \eqref{eq:introtail} therefore corresponds to an exponential cost $\exp\{-1/(3\delta^2)+o(\delta^{-2})\}$ as $\delta\downarrow0$. This is the small-time side of the coalescence problem; the $t^{-2/3}$ large-time tail is instead related to the behavior of the height covariance at small spatial separations.

Here is a heuristic explanation of the constant $1/3$. For an Airy-type variational problem with a parabolic penalty, the maximizing endpoint $T$ has a cubic-exponential displacement cost on the scale
\begin{equation}\label{eq:introendpointheuristic}
 \bP(T>x)\approx \exp\{-\tfrac43 x^3\}.
\end{equation}
The exponent can be understood from the $3/2$ right-tail exponent for one-point KPZ fluctuations: overcoming a parabolic loss of order $x^2$ requires a fluctuation of that size, at cost $\exp\{-\tfrac43(x^2)^{3/2}\}$. This reasoning also suggests that subquadratic terminal rewards should not change the leading cubic exponent, although such universality is not asserted here as a theorem.

For geodesics started at $0$ and $x$ to coalesce in unit time, imagine that their endpoints meet at a point $y$ between them. If one treats the two rare transverse displacements as though their leading exponential costs add, the resulting estimate is
\begin{equation}\label{eq:introheuristicintegral}
 \int_0^x\exp\!\left\{-\frac43\bigl(y^3+(x-y)^3\bigr)\right\}\,dy.
\end{equation}
The minimum of $y^3+(x-y)^3$ occurs at $y=x/2$, giving $\frac43\cdot\frac{x^3}{4}=\frac{x^3}{3}$. This recovers the exponent in \eqref{eq:introtail}. The two displacements are emphatically not independent, so \eqref{eq:introheuristicintegral} is only a geometric heuristic; the rigorous conclusion above comes from the covariance estimates and our identity.

There is a complementary large-time question. If $V(x)=\bE[(\fh(x)-\fh(0))^2]$, then \eqref{eq:introvariance} gives $V'(x)=2p(x)$ almost everywhere. Formally combining this with scaling yields
\begin{equation}\label{eq:introlargetime}
 \bP(\tau_{0,1}>t)=1-\frac12 V'(t^{-2/3}).
\end{equation}
For example, a differentiable second-order expansion
\begin{equation}\label{eq:introsmallx}
 V(x)=2x-cx^2+o(x^2),\qquad V'(x)=2-2cx+o(x),
\end{equation}
would imply $\bP(\tau_{0,1}>t)=ct^{-2/3}+o(t^{-2/3})$. This is a conditional observation, not a claim that the required second-order expansion has been proved. There is a point-process interpretation of its coefficient. Write $N(0,x)=\#(\mathcal P\cap(0,x))$. If the lower point field has finite intensity $\rho$ and
\begin{equation}\label{eq:intropaircondition}
\bE\big[N(0,x)(N(0,x)-1)\big]=o(x)\qquad(x\downarrow0),
\end{equation}
then stationarity and the inequalities
\[
0\le\rho x-\bP(N(0,x)\ge1)\le\tfrac12\bE[N(0,x)(N(0,x)-1)]
\]
show that $p(x)=1-\rho x+o(x)$. Consequently \eqref{eq:introsmallx} holds with $c=\rho$, and the temporal coalescence tail is asymptotic to $\rho t^{-2/3}$. Condition \eqref{eq:intropaircondition}, which would follow for instance from a locally bounded second factorial correlation density, is not established here. In particular, the temporal coalescence time in this discussion should not be confused with the spatial distance $D$ to the next lower point, whose mean is finite by \eqref{eq:intromeandistance}.

Section~2 gives the precise statements of the two main results. The remainder of the paper is organized as follows. Section~\ref{sec:KPZ integration by parts proof} proves the directed-landscape integration-by-parts identity, and Section~\ref{sec:covariance proof} derives the covariance--coalescence formula. The appendices establish the Brownian polymer and last-passage identities, justify convergence of the relevant geodesics and coalescence times, and collect the required moment estimates.

\section{Setting and main results}

We use the directed landscape, the point-to-line height $\fh$, the rightmost optimizing geodesics $\gamma(x)$, and their coalescence times $\tau_{x,y}$ as defined in the introduction. Throughout, $\bP$ and $\bE$ denote probability and expectation over the random environment. The precise choice of rightmost optimizers is used in the proofs below. Our results are stated at time one; other times follow from KPZ scaling.

\subsection{Covariance and coalescence}

For flat terminal reward, retain the notation $C$ and $p$ from \eqref{eq:introcov}. The required height moments are recorded in Lemma~\ref{lem:landscape height tails}.

\begin{theorem}[Covariance--coalescence identity]
\label{thm:covariance coalescence}
For flat terminal reward $\fh_0\equiv0$ and $0\le a<b$,
\begin{equation*}
C(a)-C(b)=\int_a^b p(x)\,dx.
\end{equation*}
In particular, $C$ is absolutely continuous on $[0,\infty)$, and
\[
C'(x)=-p(x)
\quad\text{for Lebesgue-almost every }x>0.
\]
\end{theorem}
The proof is given in Section~\ref{sec:covariance proof}.

\subsection{Flat short-distance coefficient}
\label{subsec:flatcurvature}

The following gives the short-distance curvature and the corresponding
intensity. Its analytic input is the integrated Airy$_1$ remainder estimate
proved in Appendix~\ref{app:flat-curvature}.

\begin{proposition}[Flat lower-point intensity]
\label{prop:flat-intensity}
For flat initial data at time one,
\[
\bE[(\fh(x)-\fh(0))^2]
=2x+\frac43\bE\fh(0)\,x^2+o(x^2),\qquad x\downarrow0.
\]
Consequently, the right intensity of the lower-point field, defined by
$\rho=\lim_{x\downarrow0}x^{-1}\bP(\mathcal P\cap(0,x)\neq\varnothing)$,
exists and is
\[
\rho=-\frac43\bE\fh(0)
=-\frac{2^{4/3}}3\bE\chi_{\rm GOE}.
\]
\end{proposition}
\begin{proof}
The integrated Airy$_1$ expansion in Lemma~\ref{lem:flatlayer} yields
\[
\bE[(\fh(x)-\fh(0))_+^2]
=x+2\left.\partial_t\bE\fh(t,0)\right|_{t=1}x^2+o(x^2).
\]
Flat spatial stationarity and reflection invariance imply that the positive
and negative parts of the increment have equal second moments.  Since
$\bE\fh(t,0)=t^{1/3}\bE\fh(1,0)$, we obtain the first expansion.

By Theorem~\ref{thm:covariance coalescence},
$\int_0^x(1-p(y))\,dy=x-\frac12\bE[(\fh(x)-\fh(0))^2]
=-\frac23\bE\fh(0)\,x^2+o(x^2)$.
As $1-p$ is nondecreasing, monotone difference quotients show that
$(1-p(x))/x\to-\frac43\bE\fh(0)$.
Finally $\fh(0)\stackrel{\rm d}{=}2^{1/3}\mathcal A_1(0)
\stackrel{\rm d}{=}2^{-2/3}\chi_{\rm GOE}$.
\end{proof}

\subsection{Integration by parts for the height function}

For functions $g$ and $v$, write
\[
\jbrac{g,v}=\int_\R g(x)v(x)\,dx.
\]
Let $\fu=\partial_x\fh$ in the distributional sense, so that
\[
\jbrac{g,\fu}=\jbrac{-g',\fh},\qquad g\in C_c^1(\R).
\]
Throughout the paper, a cylindrical functional and its variational derivative are written as
\begin{equation}
\begin{split}
F(v)&=G\bigl(\jbrac{g_1,v},\ldots,\jbrac{g_n,v}\bigr),\\
D_xF(v)&=\sum_{j=1}^n
\partial_jG\bigl(\jbrac{g_1,v},\ldots,\jbrac{g_n,v}\bigr)g_j(x),
\end{split}
\label{eq:cylindrical functional}
\end{equation}
where $g_1,\ldots,g_n\in C_c(\R)$ and $G\in C_b^1(\R^n)$.

\begin{theorem}[Integration by parts for the KPZ fixed point]
\label{thm:KPZ fixed point integration by parts}
Let $\fh_0\in C_c^1(\R)$, and let $F$ be a cylindrical functional as in \eqref{eq:cylindrical functional}. Then, for every $g_0\in C_c^1(\R)$,
\begin{align*}
\bE\left[F(\fh)\jbrac{g_0,\fu}\right]
&=\iint_{\R^2}g_0(x_0)
\bE\left[
D_{x_1}F(\fh)\sgn(x_1-x_0)
\ind{\tau_{x_0,x_1}\le1}
\right]\,dx_0\,dx_1\\
&\quad+\int_\R g_0(x_0)
\bE\left[F(\fh)\fh_0'(\gamma_1(x_0))\right]\,dx_0.
\end{align*}
\end{theorem}
The proof is given in Section~\ref{sec:KPZ integration by parts proof}.

\section{Proofs of the directed landscape identities}

\subsection{Integration by parts}
\label{sec:KPZ integration by parts proof}

The argument uses the corresponding integration-by-parts identity for Brownian LPP and the convergence of rescaled Brownian LPP to the directed landscape. Fix $\fh_0\in C_c^1(\R)$, and denote the rescaled Brownian LPP height, geodesics, and coalescence times by $\fh^N$, $\gamma^N$, and $\tau^N$, respectively; their precise construction is given in Appendix~\ref{sec:rescaled Brownian LPP}. Lemma~\ref{lem:rescaled Brownian integration by parts} gives, for every cylindrical functional $F$ as in \eqref{eq:cylindrical functional} and every $g_0\in C_c^1(\R)$,
\begin{align}
\bE\left[F(\fh^N)\jbrac{-g_0',\fh^N}\right]
&=\iint_{\R^2}g_0(x_0)
\bE\left[
D_{x_1}F(\fh^N)\sgn(x_1-x_0)
\ind{\tau^N_{x_0,x_1}\le1}
\right]\,dx_0\,dx_1\notag\\
&\quad+\int_\R g_0(x_0)
\bE\left[F(\fh^N)\fh_0'(\gamma_1^N(x_0))\right]\,dx_0.
\label{eq:rescaled Brownian integration by parts}
\end{align}
Lemma~\ref{lem:landscape height path and coalescence convergence} shows that there is a common coupling with the directed landscape under which, for Lebesgue-almost every $(x,y)\in\R^2$,
\begin{equation}
\left(\fh^N,\gamma_1^N(x), \tau^N_{x,y} \right)
\longrightarrow
\left(\fh,\gamma_1(x),\tau_{x,y} \right)
\quad\text{in probability},
\label{eq:joint height endpoint and coalescence convergence}
\end{equation}
where the first coordinate is equipped with the topology of locally uniform convergence.

\begin{proof}[Proof of Theorem~\ref{thm:KPZ fixed point integration by parts}]
Set
\[
X(v)=\jbrac{-g_0',v}.
\]
Choose $\eta\in C_c^\infty(\R)$ with $0\le\eta\le1$, $\eta=1$ on $[-1,1]$, and $\eta=0$ outside $[-2,2]$. For $R\ge1$, define
\[
\eta_R(z)=\eta(z/R),
\qquad F_R(v)=F(v)\eta_R(X(v)).
\]
This is again a cylindrical functional, and
\[
D_xF_R(v)=\eta_R(X(v))D_xF(v)
-F(v)\eta_R'(X(v))g_0'(x).
\]
Applying \eqref{eq:rescaled Brownian integration by parts} to $F_R$ gives
\begin{align*}
\bE\left[F_R(\fh^N)X(\fh^N)\right]
&=\iint_{\R^2}g_0(x_0)
\bE\left[
D_{x_1}F_R(\fh^N)\sgn(x_1-x_0)
\ind{\tau^N_{x_0,x_1}\le1}
\right]\,dx_0\,dx_1\\
&\quad+\int_\R g_0(x_0)
\bE\left[F_R(\fh^N)\fh_0'(\gamma_1^N(x_0))\right]\,dx_0.
\end{align*}

Keep $R$ fixed and let $N\to\infty$. The map $v\mapsto F_R(v)X(v)$ is bounded and continuous in the topology of locally uniform convergence, so
\[
\bE\left[F_R(\fh^N)X(\fh^N)\right]
\longrightarrow\bE\left[F_R(\fh)X(\fh)\right].
\]
By \eqref{eq:joint height endpoint and coalescence convergence}, the random integrands on the right-hand side converge in probability for almost every spatial argument. They are bounded by deterministic integrable functions: all test functions have compact support, $G$ and its first derivatives are bounded, and $\fh_0'$ is bounded. Dominated convergence therefore yields
\begin{align*}
\bE\left[F_R(\fh)X(\fh)\right]
&=\iint_{\R^2}g_0(x_0)
\bE\left[
D_{x_1}F_R(\fh)\sgn(x_1-x_0)
\ind{\tau_{x_0,x_1}\le1}
\right]\,dx_0\,dx_1\\
&\quad+\int_\R g_0(x_0)
\bE\left[F_R(\fh)\fh_0'(\gamma_1(x_0))\right]\,dx_0.
\end{align*}

Now let $R\to\infty$. Lemma~\ref{lem:landscape height tails} gives
\[
\bE|X(\fh)|
\le\|g_0'\|_{L^1}\sup_{x\in\R}\bE|\fh(x)|<\infty.
\]
Since $|F_R(\fh)|\le\|G\|_\infty$, dominated convergence applies to the left-hand side. In the first term on the right-hand side, the contribution involving $\eta_R'$ is bounded in absolute value by
\[
\frac{1}{R}\|G\|_\infty\|\eta'\|_\infty
\|g_0\|_{L^1}\|g_0'\|_{L^1},
\]
and hence vanishes. The remaining terms converge by dominated convergence. Since $X(\fh)=\jbrac{g_0,\fu}$, the resulting identity is precisely the statement of the theorem.
\end{proof}

\subsection{Proof of the covariance identity}
\label{sec:covariance proof}
Throughout this subsection, we set $\fh_0\equiv0$.

\begin{lemma}[Height increments]
\label{lem:height increments}
For $x_0,x_1,x_2\in\R$,
\[
\bE\left[\fh(x_2)\bigl(\fh(x_1)-\fh(x_0)\bigr)\right]
=\int_{x_0}^{x_1}\sgn(x_2-x)\bP(\tau_{x_2,x}\le1)\,dx.
\]
\end{lemma}

\begin{proof}
Assume first that $x_0<x_1$. Let $\rho\in C_c^\infty(\R)$ be nonnegative and even, with $\int_\R\rho=1$, and set
\[
\rho_\epsilon(x)=\epsilon^{-1}\rho(x/\epsilon),
\qquad\fh_\epsilon=\rho_\epsilon*\fh,
\qquad g^\epsilon=\rho_\epsilon*\ind{[x_0,x_1]}.
\]
To apply Theorem~\ref{thm:KPZ fixed point integration by parts} to the unbounded observable $\fh_\epsilon(x_2)$, choose $\eta\in C_c^\infty(\R)$ with $0\le\eta\le1$ and $\eta=1$ near $0$, and define
\[
G_R(z)=z\eta(z/R),
\qquad F_R(v)=G_R\bigl(\jbrac{\rho_\epsilon(x_2-\cdot),v}\bigr).
\]
Theorem~\ref{thm:KPZ fixed point integration by parts} gives
\begin{align}
\bE\left[F_R(\fh)\jbrac{g^\epsilon,\fu}\right]
&=\iint_{\R^2}g^\epsilon(z_0)
\bE\left[
G_R'(\fh_\epsilon(x_2))\rho_\epsilon(x_2-z_1)
\sgn(z_1-z_0)\ind{\tau_{z_1,z_0}\le1}
\right]\,dz_0\,dz_1.
\label{eq:height increment ibp}
\end{align}

We first let $\epsilon\downarrow0$ with $R$ fixed. The deterministic point $x_2$ is almost surely not a jump point of the piecewise constant endpoint map $x\mapsto\gamma_1(x)$; see \cite[Lemma 3.6 and Proposition 3.10]{rahman2025infinite}, applied in our point-to-line convention by reversing space and time. Thus, almost surely, there is a $\delta>0$ such that $\gamma_1(z)=\gamma_1(x_2)$ whenever $|z-x_2|<\delta$. On this event,
\[
\ind{\tau_{z,z_0}\le1}=\ind{\tau_{x_2,z_0}\le1},
\qquad |z-x_2|<\delta,
\]
for every $z_0\in\R$. Set
\[
J_\epsilon(z_0)=\int_\R\rho_\epsilon(x_2-z_1)
\sgn(z_1-z_0)\ind{\tau_{z_1,z_0}\le1}\,dz_1.
\]
For every $z_0\ne x_2$, local constancy gives, almost surely,
\[
J_\epsilon(z_0)\longrightarrow
\sgn(x_2-z_0)\ind{\tau_{x_2,z_0}\le1},
\qquad |J_\epsilon(z_0)|\le1.
\]
By Fubini, the right-hand side of \eqref{eq:height increment ibp} equals
\[
\bE\left[G_R'(\fh_\epsilon(x_2))
\int_\R g^\epsilon(z_0)J_\epsilon(z_0)\,dz_0\right].
\]
Replacing $g^\epsilon$ by $\ind{[x_0,x_1]}$ changes this expression by at most
\[
\|G_R'\|_\infty
\|g^\epsilon-\ind{[x_0,x_1]}\|_{L^1}\longrightarrow0.
\]
Continuity of $\fh$ gives $G_R'(\fh_\epsilon(x_2))\to G_R'(\fh(x_2))$ almost surely. Dominated convergence therefore shows that the right-hand side of \eqref{eq:height increment ibp} tends to
\[
\int_{x_0}^{x_1}\sgn(x_2-z_0)
\bE\left[G_R'(\fh(x_2))\ind{\tau_{x_2,z_0}\le1}\right]\,dz_0.
\]

For the left-hand side, note that
\[
\jbrac{g^\epsilon,\fu}
=\jbrac{-(g^\epsilon)',\fh}
=\fh_\epsilon(x_1)-\fh_\epsilon(x_0).
\]
By continuity, its random integrand converges almost surely to
\[
G_R(\fh(x_2))\bigl(\fh(x_1)-\fh(x_0)\bigr).
\]
Jensen's inequality and Lemma~\ref{lem:landscape height tails} give
\[
\sup_{0<\epsilon\le1}\sup_{x\in\R}
\bE|\fh_\epsilon(x)|^4
\le\sup_{x\in\R}\bE|\fh(x)|^4<\infty.
\]
Since $|G_R(z)|\le|z|$, Cauchy--Schwarz implies a uniform second-moment bound on
\[
G_R(\fh_\epsilon(x_2))
\bigl(\fh_\epsilon(x_1)-\fh_\epsilon(x_0)\bigr),
\qquad R\ge1,\quad 0<\epsilon\le1.
\]
These random variables are therefore uniformly integrable. Taking $\epsilon\downarrow0$ in \eqref{eq:height increment ibp}, we obtain
\begin{align*}
&\bE\left[G_R(\fh(x_2))\bigl(\fh(x_1)-\fh(x_0)\bigr)\right]\\
&\qquad=\int_{x_0}^{x_1}\sgn(x_2-z_0)
\bE\left[G_R'(\fh(x_2))\ind{\tau_{x_2,z_0}\le1}\right]\,dz_0.
\end{align*}
Finally, $G_R(z)\to z$, $G_R'(z)\to1$, and $\sup_{R\ge1}\|G_R'\|_\infty<\infty$. The left-hand side is dominated by the integrable random variable
\[
|\fh(x_2)|\bigl(|\fh(x_1)|+|\fh(x_0)|\bigr).
\]
Dominated convergence on both sides proves the identity for $x_0<x_1$. The case $x_0=x_1$ is immediate, and the case $x_0>x_1$ follows by interchanging $x_0$ and $x_1$.
\end{proof}

\begin{proof}[Proof of Theorem~\ref{thm:covariance coalescence}]
Under flat terminal reward, spatial stationarity gives $\bE\fh(a)=\bE\fh(b)$. For $0\le a<b$, Lemma~\ref{lem:height increments}, with $x_2=0$, $x_0=a$, and $x_1=b$, therefore yields
\[
C(b)-C(a)
=\bE\left[\fh(0)\bigl(\fh(b)-\fh(a)\bigr)\right]
=-\int_a^b p(x)\,dx.
\]
Since $0\le p\le1$, absolute continuity and the almost-everywhere derivative formula follow.
\end{proof}

\appendix

\section{Brownian LPP integration by parts}

This appendix derives the Brownian LPP integration-by-parts identity used in Lemma~\ref{lem:rescaled Brownian integration by parts}. We first work with a fixed number of layers $N$ and a general terminal reward $h_0$. The KPZ scaling and the choice of reward approximating $\fh_0$ are specified in Section~\ref{sec:rescaled Brownian LPP}.

\subsection{O'Connell--Yor polymer}

Let $\{B^{(k)}\}_{k\ge1}$ be independent standard two-sided Brownian motions, with joint law and expectation denoted by $\bP$ and $\bE$. Let $w_1,\ldots,w_N$ be independent exponential random variables with rate one; expectation over these variables is denoted by $\E_w$.

For $r\in\R$, set
\[
T_0^w(r)=r,\qquad T_k^w(r)=r+\sum_{j=1}^k w_j,
\quad 1\le k\le N.
\]
Write
\[
\Delta_N(r,s)=\{(t_0,\ldots,t_N):r=t_0\le\cdots\le t_N=s\},
\qquad \Delta_N(r)=\bigcup_{s\ge r}\Delta_N(r,s), \qquad \Delta_N = \bigcup_{r} \Delta_N(r).
\]
Fix a terminal reward $h_0\in C_c^1(\R)$ and inverse temperature $0<\beta<\infty$. For $\bt\in\Delta_N(r)$, define the path weight, including the terminal reward, by
\begin{equation*}
H(\bt)
=\sum_{k=1}^N\bigl(B^{(k)}(t_k)-B^{(k)}(t_{k-1})\bigr)
-N-(t_N-t_0)+h_0(t_N).
\end{equation*}
The dependence of $H$ on $N$ and $h_0$ is implicit. Define the point-to-line partition function and free energy by
\begin{align*}
Z^\beta(r)
&=\E_w \sqbrac{
e^{\beta \pare{H(T^w(r))}}
},
\\
h^\beta(r)&=\beta^{-1}\log Z^\beta(r).
\end{align*}
For a bounded measurable path observable $\Psi$, define the quenched polymer measure starting at $r$ by
\[
\Epoly_r\Psi(T)
=\frac{\E_w\left[
 e^{\beta[H(T^w(r))]}\Psi(T^w(r))\right]}
 {Z^\beta(r)}.
\]
Write $\Ppoly_r$ for the corresponding quenched probability. The notation $\Epoly_{r_0,\ldots,r_m}$ and $\Ppoly_{r_0,\ldots,r_m}$ always refers to conditionally independent paths in the same environment. Here $T$ is the path under the polymer measure; under a product measure, $T^j$ starts at $r_j$, and $T^{j'}$ denotes an independent replica at the same start.

For two paths $\bs,\bt \in \Delta_N$, define their first intersection level by
\begin{equation*}
\kappa(\bs,\bt)
=\inf\{k\in\{1,\ldots,N\}:
[s_{k-1},s_k]\cap[t_{k-1},t_k]\ne\varnothing\},
\qquad\inf\varnothing=\infty.
\end{equation*}
For sampled paths $T^0,T^1$ starting at $r_0,r_1$, write
\[
\kappa_{r_0,r_1}^\beta=\kappa(T^0,T^1).
\]
At finite temperature this records the first intersection; the two sampled paths need not coincide afterwards.

\subsection{Brownian last passage and optimizing paths}

For integer layers $1\le m\le n$ and $r\le s$, define the normalized point-to-point last-passage value by
\begin{equation}
L(r,m;s,n)
=\sup_{r=t_{m-1}\le t_m\le\cdots\le t_n=s}
\left\{\sum_{k=m}^n\bigl(B^{(k)}(t_k)-B^{(k)}(t_{k-1})\bigr)\right\}
-(n-m+1)-(s-r).
\label{eq:Brownian point to point value}
\end{equation}
The point-to-line value with reward $h_0$ is
\begin{equation}
L_N(r)=\sup_{s\ge r}\{L(r,1;s,N)+h_0(s)\}
=\sup_{\bt\in\Delta_N(r)}H(\bt).
\label{eq:Brownian point to line value}
\end{equation}
When $N$ and $h_0$ are fixed, we abbreviate $L_N$ to $L$.

Let $\Gamma(r)=(\Gamma_k(r))_{k=0}^N$ be the rightmost optimizing path for $L(r)$. Thus $\Gamma_0(r)=r$, and $\Gamma_N(r)$ is its maximizing terminal location. The dependence on $N$ and $h_0$ is implicit. Define
\[
\kappa_{r_0,r_1}^\infty=\kappa\bigl(\Gamma(r_0),\Gamma(r_1)\bigr).
\]
We call a starting point \emph{non-exceptional} if its optimizing path is unique. A finite first intersection is \emph{non-degenerate} if the two path intervals at that level have an intersection of positive length.

\begin{lemma}[Brownian LPP optimizing paths]
Fix a deterministic terminal reward $h_0\in C_c^1(\R)$. For each $N\ge1$, let $\Gamma(r)$ denote the rightmost optimizing path for the Brownian point-to-line last-passage value $L_N(r)$.
\begin{enumerate}
    \item For sufficiently large $N$, depending only on $h_0$, the optimizing path from each fixed deterministic starting point $r$ is almost surely unique. Moreover, almost surely, the set of starting points admitting more than one optimizing path has Lebesgue measure zero.
    \item For every fixed $N\ge1$ and every deterministic pair $r\neq s$, the first intersection of $\Gamma(r)$ and $\Gamma(s)$, when finite, is almost surely non-degenerate.
\end{enumerate}
\end{lemma}

\begin{proof}
Part~1 follows from \cite[Lemma~4.6(2)]{hammond2019patchwork}.

We now prove part~2. By symmetry, fix deterministic $r<s$, and write $a_j=\Gamma_j(r)$ and $b_j=\Gamma_j(s)$. For $1\le j<q\le N+1$, define
\[
F_j^q(t)=
\begin{cases}
B^{(j)}(t)-B^{(q)}(t), & q\le N,\\
B^{(j)}(t)-t+h_0(t), & q=N+1.
\end{cases}
\]
For distinct $i,j<q\le N$, the pair $(F_i^q,F_j^q)$ is a two-dimensional Brownian motion with covariance matrix per unit time
\[
\begin{pmatrix}2&1\\1&2\end{pmatrix}.
\]
For $q=N+1$, the pair consists of two independent Brownian motions with the common deterministic drift $h_0(t)-t$. Representing each pair in terms of two independent standard Brownian motions, \cite[Theorem~1.1]{cator2012local} implies that almost surely $F_i^q$ and $F_j^q$ have no common local maximum. For the terminal case, we apply the theorem after cutting off the drift outside bounded intervals. Moreover, almost surely $F_j^q$ has a local maximum not equal to the deterministic point $s$. We work on the probability-one event where these properties holds for every choice of indices.

Suppose, for a contradiction, that the first intersection occurs on row $k$ and consists of a single point $z$. Exchanging optimal paths after the intersection preserves optimality, so $a_j=b_j$ for all $j\ge k$. Consequently,
\[
a_k=b_{k-1}=b_k=z.
\]
Let
\[
q=\min\{j>k:a_j>z\},
\]
with $q=N+1$ if this set is empty. Then
\[
a_k=\cdots=a_{q-1}=b_k=\cdots=b_{q-1}=z,
\]
and, if $q\le N$, $a_q=b_q>z$. We also have $a_{k-1}<z$. Indeed, for $k=1$ this follows from $a_0=r<s=z$, while for $k\ge2$ disjointness on row $k-1$ gives $a_{k-1}<b_{k-2}\le z$. For either sign of sufficiently small $\epsilon$, replacing $a_k,\ldots,a_{q-1}$ by $z+\epsilon$ therefore gives an admissible path. Its weight changes by
\[
F_k^q(z+\epsilon)-F_k^q(z).
\]
Optimality implies that $z$ is a local maximum of $F_k^q$.

If $z=s$, this contradicts our choice of the probability-one event. Otherwise, $z>s$. Define
\[
p=\min\{j\ge1:b_j=z\}.
\]
Since $b_0=s<z$ and $b_{k-1}=z$, we have $p<k$, $b_{p-1}<z$, and
\[
b_p=\cdots=b_{q-1}=z.
\]
Thus, for either sign of sufficiently small $\epsilon$, replacing $b_p,\ldots,b_{q-1}$ by $z+\epsilon$ also gives an admissible path, whose weight changes by
\[
F_p^q(z+\epsilon)-F_p^q(z).
\]
Optimality makes $z$ a local maximum of $F_p^q$ as well. This contradicts the absence of a common local maximum of $F_k^q$ and $F_p^q$.
\end{proof}

The following convergence from the OY polymer to Brownian LPP follows from the Laplace principle.
\begin{lemma}[OY to LPP]\label{lem:polymer to Brownian LPP}
Almost surely in the Brownian environment, the following statements hold as $\beta\to\infty$.
\begin{enumerate}
\item For any compact set $K$,
\[
\sup_{r\in K} |h^\beta(r) - L(r)| \longrightarrow 0.
\]
\item For any non-exceptional $r$ and every $\delta>0$,
\[
\Ppoly_r\left(\max_{0\le k\le N}|T_k-\Gamma_k(r)|>\delta\right)
\longrightarrow0.
\]

\item For any non-exceptional $r_0,r_1$ whose first intersection, when finite, is non-degenerate,
\[
\Ppoly_{r_0,r_1}\left(\kappa^\beta_{r_0,r_1}\ne\kappa^\infty_{r_0,r_1}\right)
\longrightarrow0.
\]
\end{enumerate}
\end{lemma}

\subsection{Gaussian integration by parts}

\begin{lemma}[Regularization of the polymer environment]
\label{lem:polymer environment regularization}
Fix $N<\infty$, $0<\beta<\infty$, and $h_0\in C_c^1(\R)$. Let $\rho\in C_c^\infty(\R)$ be even and nonnegative, with $\int_\R\rho=1$, and set, for $0 \le \epsilon \le 1$,
\[
\rho_\epsilon(t)=\epsilon^{-1}\rho(t/\epsilon),
\qquad B_\epsilon^{(k)}=\rho_\epsilon*B^{(k)}.
\]
Define $H_\epsilon$, $Z_\epsilon^\beta$, $h_\epsilon^\beta$, and $\Ppoly_{r;\epsilon}$ by replacing every $B^{(k)}$ by $B_\epsilon^{(k)}$ in the weight, partition function, free energy, and polymer measure. Then $h_\epsilon^\beta\in C^1(\R)$ almost surely, and
\begin{equation}\label{eq:regularized spatial derivative}
(h_\epsilon^\beta)'(r)
=\Epoly_{r;\epsilon}\left[
\sum_{k=1}^N\bigl((B_\epsilon^{(k)})'(T_k)
                 -(B_\epsilon^{(k)})'(T_{k-1})\bigr)
+h_0'(T_N)\right].
\end{equation}
For every compact interval $K$ and every $1\le p<\infty$,
\begin{align*}
\sup_{r\in K}|h_\epsilon^\beta(r)-h^\beta(r)|&\longrightarrow0
\quad\text{almost surely and in }L^p(\bP),
\\
\sup_{r\in K}\|\Ppoly_{r;\epsilon}-\Ppoly_r\|_{\mathrm{TV}}
&\longrightarrow0\quad\text{almost surely}.
\end{align*}
\end{lemma}

\begin{proof}
We first prove convergence. Set $B_0^{(k)}=B^{(k)}$, with the same convention for the associated quantities. Fix a compact interval $K\subset[-R_K,R_K]$, and choose $R_\rho>0$ such that $\operatorname{supp}\rho\subset[-R_\rho,R_\rho]$. Define
\[
M=\max_{1\le k\le N}\sup_{t\in\R}\bigl(|B^{(k)}(t)|-\frac{1}{4N}|t|\bigr).
\]
The standard exponential-martingale estimate for Brownian motion, together with a union bound, gives
\[
\bP(M>u)\le4N e^{-\frac{u}{2N}},\qquad u\ge0.
\]
Thus $M$ is almost surely finite and has moments of every order. Since $\rho$ is nonnegative and has integral one,
\[
|B_\epsilon^{(k)}(t)|\le M+\frac{1}{4N}(|t|+R_\rho),
\qquad 0\le\epsilon\le1.
\]
Let $\bt\in\Delta_N(r)$ with $r\in K$. Every coordinate of $\bt$ has absolute value at most $R_K+t_N-r$, so
\[
\left|\sum_{k=1}^N\bigl(B_\epsilon^{(k)}(t_k)-B_\epsilon^{(k)}(t_{k-1})\bigr)\right|
\le2NM+\frac12(R_K+R_\rho+t_N-r).
\]
Consequently, the nonnegative random variable
\[
A_K=2NM+N+\|h_0\|_\infty+\frac12(R_K+R_\rho)
\]
has an exponential tail and satisfies
\begin{equation}\label{eq:uniform polymer weight bounds}
-A_K-\frac32(t_N-r)\le H_\epsilon(\bt)\le A_K-\frac12(t_N-r),
\qquad r\in K,\quad \bt\in\Delta_N(r),\quad 0\le\epsilon\le1.
\end{equation}

Almost surely, $B_\epsilon^{(k)}\to B^{(k)}$ uniformly on compact intervals. Hence, for each fixed $w$,
\[
\sup_{r\in K}|H_\epsilon(T^w(r))-H(T^w(r))|\longrightarrow0.
\]
Writing $S=T_N^w(r)-r=\sum_{j=1}^Nw_j$, \eqref{eq:uniform polymer weight bounds} gives
\[
\sup_{0\le \epsilon \le 1}\sup_{r\in K}\left|e^{\beta H_\epsilon(T^w(r))}\right|
\le e^{\beta A_K}e^{-\beta S/2}.
\]
 Dominated convergence therefore yields
\[
D_\epsilon:=\E_w\left[\sup_{r\in K}\left|e^{\beta H_\epsilon(T^w(r))}-e^{\beta H(T^w(r))}\right|\right]\longrightarrow0.
\]
Since $\sup_{r\in K}|Z_\epsilon^\beta(r)-Z^\beta(r)|\le D_\epsilon$, the partition functions converge uniformly on $K$. 

Since $\E_w[T_N^w(r)-r]=N$, Jensen's inequality and the upper bound above give
\[
e^{-\beta(A_K+3N/2)}\le Z_\epsilon^\beta(r)\le e^{\beta A_K}.
\]
Hence,
\[
\sup_{0\le\epsilon\le1}\sup_{r\in K}|h_\epsilon^\beta(r)|
\le A_K+\frac32N.
\]
The uniform convergence of partition functions and this uniform bound give uniform convergence of $h_\epsilon^\beta$ to $h^\beta$ on $K$. And the moment bounds for $A_K$ upgrade this convergence to $L^p(\bP)$ for every $1\le p<\infty$.

For each $r$, represent both polymer measures by their densities relative to the common reference law of $w$. Comparing the normalized densities gives
\[
\sup_{r\in K}\|\Ppoly_{r;\epsilon}-\Ppoly_r\|_{\mathrm{TV}}
\le\frac{2D_\epsilon}{\inf_{r\in K}Z^\beta(r)}
\longrightarrow0
\]
almost surely.

Next, we prove \eqref{eq:regularized spatial derivative}. Differentiating the convolution and using the linear-growth bound on the Brownian paths gives
\[
|(B_\epsilon^{(k)})'(t)|\le C(1+M+|t|).
\]
Every coordinate $T_k^w(r)$ has derivative one with respect to $r$, whereas $T_N^w(r)-r=S$ is independent of $r$. Therefore,
\[
\partial_rH_\epsilon(T^w(r))
=\sum_{k=1}^N\bigl((B_\epsilon^{(k)})'(T_k^w(r))-(B_\epsilon^{(k)})'(T_{k-1}^w(r))\bigr)+h_0'(T_N^w(r)).
\]
Uniformly over $r\in K$, its absolute value is bounded by $C_{K}(1+M+S)$. The weight estimate \eqref{eq:uniform polymer weight bounds}, therefore justifies differentiation under $\E_w$ and continuity of the resulting derivative. Thus $h_\epsilon^\beta\in C^1(\R)$, and differentiating $\beta^{-1}\log Z_\epsilon^\beta$ gives \eqref{eq:regularized spatial derivative}.
\end{proof}

\begin{lemma}[Brownian LPP integration by parts]
\label{lem:Brownian LPP integration by parts}
Let $h_0\in C_c^1(\R)$, $f_1,\ldots,f_n\in C_c(\R)$, and $G\in C_b^1(\R^n)$. Define
\[
F(v)=G\bigl(\jbrac{f_1,v},\ldots,\jbrac{f_n,v}\bigr),
\]
and use the derivative $D_rF(v)$ defined in \eqref{eq:cylindrical functional}. Then, for every $f_0\in C_c^1(\R)$,
\begin{align}
\bE[F(L)\jbrac{-f_0',L}]
&=\frac{1}{2}\iint_{\R^2}f_0(r_0)
\bE\left[D_{r_1}F(L)
\sgn(r_1-r_0)\ind{\kappa_{r_0,r_1}^\infty\le N}\right]
 dr_0\,dr_1\notag\\
&\quad+\int_\R f_0(r_0)
\bE\left[F(L)h_0'(\Gamma_N(r_0))\right]dr_0.
\label{eq:Brownian integration by parts}
\end{align}
\end{lemma}

\begin{proof}
Fix $\beta<\infty$ and use the regularization of Lemma~\ref{lem:polymer environment regularization}. We first perform Gaussian integration by parts at fixed $\epsilon>0$.

Let $\cD_k$ denote the Malliavin derivative with respect to the white noise $dB^{(k)}$, taking values in $L^2(\R)$, and use the $L^2(\R)$ inner product below. Set
\[
v_{k,\epsilon}(t;\bt)
=\rho_\epsilon(t-t_k)-\rho_\epsilon(t-t_{k-1}),\qquad
W_\epsilon(\bt)=\frac{e^{\beta H_\epsilon(\bt)}}{Z_\epsilon^\beta(t_0)}.
\]
Applying spatial integration by parts to \eqref{eq:regularized spatial derivative} and then Gaussian integration by parts gives
\begin{align}
\bE[F(h_\epsilon^\beta)\jbrac{-f_0',h_\epsilon^\beta}]
&=\int_\R f_0(r_0)\bE\E_{w^0}\left[
\sum_{k=1}^N
\left\langle
\cD_k\left(W_\epsilon(T^{w^0}(r_0))F(h_\epsilon^\beta)\right),
v_{k,\epsilon}(\cdot;T^{w^0}(r_0))
\right\rangle\right]dr_0\notag\\
&\quad+\int_\R f_0(r_0)\bE\left[
F(h_\epsilon^\beta)\Epoly_{r_0;\epsilon}h_0'(T_N^0)
\right]dr_0.
\label{eq:Gaussian integration by parts}
\end{align}
Set
\[
q_\epsilon=\rho_\epsilon*\rho_\epsilon,
\]
and define
\begin{align}\label{eq:regularized intersection kernel}
I^\epsilon_k(\bs,\bt)
&=\left\langle
\rho_\epsilon*\ind{[s_{k-1},s_k]},
v_{k,\epsilon}(\cdot;\bt)\right\rangle \nonumber\\
&=(q_\epsilon*\ind{[s_{k-1},s_k]})(t_k)
 -(q_\epsilon*\ind{[s_{k-1},s_k]})(t_{k-1}).
\end{align}
Evenness of $q_\epsilon$ gives
\begin{equation*}
I^\epsilon_k(\bs,\bt)
=-I^\epsilon_k(\bt,\bs),
\qquad |I^\epsilon_k|\le1.
\end{equation*}
The quotient rule and the chain rule show that the expectation $\E_{w^0}$ in \eqref{eq:Gaussian integration by parts} equals
\begin{align*}
&\beta F(h_\epsilon^\beta)
\sum_{k=1}^N\Epoly_{r_0,r_0;\epsilon}\left[
I^\epsilon_k(T^0,T^0)-I^\epsilon_k(T^{0'},T^0)\right]\\
&\qquad+\int_\R D_{r_1}F(h_\epsilon^\beta)
\sum_{k=1}^N\Epoly_{r_0,r_1;\epsilon}
I^\epsilon_k(T^1,T^0)\,dr_1.
\end{align*}

The first line vanishes exactly: $I^\epsilon_k(T^0,T^0)=0$, and $T^0,T^{0'}$ are exchangeable under $\Ppoly_{r_0,r_0;\epsilon}$. Substitution gives
\begin{align*}
\bE[F(h_\epsilon^\beta)\jbrac{-f_0',h_\epsilon^\beta}]
&=\iint_{\R^2}f_0(r_0)\bE\left[
D_{r_1}F(h_\epsilon^\beta)
\sum_{k=1}^N\Epoly_{r_0,r_1;\epsilon}
I^\epsilon_k(T^1,T^0)\right]dr_0\,dr_1\\
&\quad+\int_\R f_0(r_0)\bE\left[
F(h_\epsilon^\beta)\Epoly_{r_0;\epsilon}h_0'(T_N^0)
\right]dr_0.
\end{align*}
The two integrands on the right-hand side are bounded, uniformly in $\epsilon$, by
\[
N|f_0(r_0)| \sum_{j=1}^n\|\partial_jG\|_\infty|f_j(r_1)|, \qquad
|f_0(r_0)|\|G\|_\infty\|h_0'\|_\infty.
\]

Next, let $\epsilon\downarrow0$ with $\beta$ fixed. Use the symmetric endpoint convention
\[
\ind{[u,v]}^*(z)
=\ind{u<z<v}+\frac12\ind{z=u}+\frac12\ind{z=v},
\qquad u<v,
\]
and define
\[
I_k(\bs,\bt)
=\ind{[s_{k-1},s_k]}^*(t_k)
 -\ind{[s_{k-1},s_k]}^*(t_{k-1}).
\]
Since $q_\epsilon$ is an even approximate identity, $I^\epsilon_k(\bs,\bt)\to I_k(\bs,\bt)$ whenever the interval $[s_{k-1}, s_k]$ has positive length. The bound $|I^\epsilon_k|\le1$ and total variation convergence of the polymer measures in Lemma~\ref{lem:polymer environment regularization} imply, almost surely,
\[
\Epoly_{r_0,r_1;\epsilon}I^\epsilon_k(T^1,T^0)
\longrightarrow\Epoly_{r_0,r_1}I_k(T^1,T^0).
\]
By antisymmetry of $I_k$ and \eqref{eq:path exchange contribution},
\[
\sum_{k=1}^N\Epoly_{r_0,r_1}I_k(T^1,T^0)
=\frac12\sgn(r_1-r_0)
\Ppoly_{r_0,r_1}(\kappa^\beta_{r_0,r_1}\le N).
\]
Applying dominated convergence, we obtain the finite-temperature identity
\begin{align}
\bE[F(h^\beta)\jbrac{-f_0',h^\beta}]
&=\frac12\iint_{\R^2}f_0(r_0)\bE\left[
D_{r_1}F(h^\beta)\sgn(r_1-r_0)
\Ppoly_{r_0,r_1}(\kappa^\beta_{r_0,r_1}\le N)
\right]dr_0\,dr_1\notag\\
&\quad+\int_\R f_0(r_0)\bE\left[
F(h^\beta)\Epoly_{r_0}h_0'(T_N^0)\right]dr_0.
\label{eq:polymer integration by parts}
\end{align}

Finally, let $\beta\to\infty$. Let $K$ be a compact set containing the supports of $f_0,f_1,\ldots,f_n$. Lemma~\ref{lem:Brownian LPP and polymer tails}, with $q=2$, gives uniform integrability of the bound
\[
|F(h^\beta)\jbrac{-f_0',h^\beta}|
\le\|G\|_\infty\|f_0'\|_{L^1}\sup_{r\in K}|h^\beta(r)|.
\]
By Lemma~\ref{lem:polymer to Brownian LPP}, $h^\beta\to L$ uniformly on $K$ almost surely. Consequently,
\[
F(h^\beta)\jbrac{-f_0',h^\beta}\longrightarrow F(L)\jbrac{-f_0',L}
\quad\text{almost surely}.
\]
Hence
\[
\lim_{\beta\to\infty}
\bE\left[F(h^\beta)\jbrac{-f_0',h^\beta}\right]
=
\bE\left[F(L)\jbrac{-f_0',L}\right].
\]
At each non-exceptional point $r_0$, convergence of the terminal point in Lemma~\ref{lem:polymer to Brownian LPP} and the fact that $h_0'\in C_c(\R)$ give
\[
\Epoly_{r_0}h_0'(T_N^0)\longrightarrow h_0'(\Gamma_N(r_0)).
\]
For any non-exceptional $r_0,r_1$ whose first intersection, if finite, is non-degenerate, Lemma~\ref{lem:polymer to Brownian LPP} gives
\begin{equation*}
\Ppoly_{r_0,r_1}(\kappa^\beta_{r_0,r_1}\le N)
\longrightarrow\ind{\kappa^\infty_{r_0,r_1}\le N}
\quad\text{almost surely},
\end{equation*}
and
\[
D_{r_1}F(h^\beta)
\longrightarrow D_{r_1}F(L).
\]
Together with the bound
\begin{equation*}
\left|D_rF(v)\right|
\le \sum_{j=1}^n\|\partial_jG\|_\infty|f_j(r)|,
\end{equation*}
these convergences allow us to apply dominated convergence to the right-hand side of \eqref{eq:polymer integration by parts}, proving \eqref{eq:Brownian integration by parts}.
\end{proof}

\begin{lemma}[Exchange at the first intersection]
Fix $N<\infty$ and $0<\beta<\infty$. Let $T^0,T^1$ be two paths under $\Ppoly_{r_0,r_1}$, and write $\kappa=\kappa_{r_0,r_1}^\beta$. Define the exchange transformation $\cS: (T^0,T^1) \mapsto (\widehat T^0, \widehat T^1)$ by
\[
\widehat T_j^i=
\begin{cases}
T_j^i,&j<\kappa,\\
T_j^{1-i},&j\ge\kappa,
\end{cases}
\qquad i=0,1,\quad j=0,\ldots,N,
\]
on $\{\kappa\le N\}$, and let $\cS$ be the identity on $\{\kappa=\infty\}$. Then $\cS$ is measure-preserving under $\Ppoly_{r_0,r_1}$.

Let $I^\epsilon_k$ be the function defined in \eqref{eq:regularized intersection kernel}. Then, for any $\epsilon>0$,
\begin{equation}
\Epoly_{r_0,r_1}\sqbrac{
\ind{k>\kappa} I_k^\epsilon(T^0,T^1)
}=0,
\label{eq:path exchange after}
\end{equation}
and
\begin{equation}
\lim_{\epsilon\downarrow0}
\Epoly_{r_0,r_1}\sqbrac{
\sum_{k=1}^N I_k^\epsilon(T^0,T^1)
}
=
-\frac12\sgn(r_1-r_0)
\Ppoly_{r_0,r_1}(\kappa\le N).
\label{eq:path exchange contribution}
\end{equation}
\end{lemma}

\begin{proof}
With the reward-inclusive weight $H$ defined above, $(T^0,T^1)$ has density
\[
\frac{\exp\{\beta[H(T^0)+H(T^1)]
                -(T_N^0-r_0)-(T_N^1-r_1)\}}
     {Z^\beta(r_0)Z^\beta(r_1)}
\prod_{i=0}^1\prod_{j=1}^N dT_j^i.
\]
On $\{\kappa\le N\}$, both $T_\kappa^0$ and $T_\kappa^1$ exceed $\max\{T_{\kappa-1}^0,T_{\kappa-1}^1\}$ almost surely, so the path exchange $\cS$ preserves $\kappa$. On each event $\{\kappa=k\}$, $\cS$ is a coordinate permutation, and the density is invariant under $\cS$. This proves the first assertion.

When $k > \kappa$, we have $(\widehat T^0_k, \widehat T^1_k) = (T^1_k, T^0_k)$, and hence $I^\epsilon_k(T^0, T^1) = I^\epsilon_k(\widehat T^1, \widehat T^0)$. We have
\[
\Epoly_{r_0,r_1}\sqbrac{
\ind{k>\kappa} I_k^\epsilon(T^0,T^1)
}
= \Epoly_{r_0,r_1}\sqbrac{
\ind{k>\kappa} I_k^\epsilon(\widehat T^1, \widehat T^0)
}
= \Epoly_{r_0,r_1}\sqbrac{
\ind{k>\kappa} I_k^\epsilon(T^1,T^0)
},
\]
where the second equality uses invariance under $\cS$.

On the other hand, $I_k^\epsilon$ is antisymmetric, that is,
\[
I^\epsilon_k(T^0, T^1) = - I^\epsilon_k(T^1, T^0).
\]
We have
\[
\Epoly_{r_0,r_1}\sqbrac{
\ind{k>\kappa} I_k^\epsilon(T^0,T^1)
}
= -\Epoly_{r_0,r_1}\sqbrac{
\ind{k>\kappa} I_k^\epsilon(T^0,T^1)
} = 0.
\]

Suppose first that $r_0<r_1$. Before the first intersection, the intervals are disjoint, so $I^\epsilon_k(T^0, T^1)$ tends to zero as $\epsilon\downarrow0$. At the first intersection, write
\[
a=T_{\kappa-1}^0,
\quad c=T_{\kappa-1}^1,
\quad b=T_\kappa^0,
\quad d=T_\kappa^1.
\]
Almost surely on $\{\kappa\le N\}$, $a<c<\min\{b,d\}$ and $b\ne d$. Consequently,
\[
\lim_{\epsilon\downarrow0}I_\kappa^\epsilon(T^0,T^1)
=-\ind{b<d},\qquad
\lim_{\epsilon\downarrow0}I_\kappa^\epsilon(\widehat T^0,\widehat T^1)
=-\ind{d<b}.
\]
The sum of these two limits is $-1$. Invariance under $\cS$ therefore makes the expected first-intersection contribution $-\tfrac12\Ppoly_{r_0,r_1}(\kappa\le N)$. Since $N$ is fixed, dominated convergence justifies all the limits. Together with \eqref{eq:path exchange after}, this proves \eqref{eq:path exchange contribution} when $r_0<r_1$. The case $r_0>r_1$ follows by exchanging the path labels. If $r_0=r_1$, the entire pair is exchangeable and antisymmetry makes the expectation zero for every $\epsilon>0$, which agrees with the convention $\sgn(0)=0$.
\end{proof}

\section{Convergence to the directed landscape}
\label{sec:directed landscape convergence}

\subsection{Rescaled Brownian LPP}
\label{sec:rescaled Brownian LPP}

Fix the terminal reward $\fh_0\in C_c^1(\R)$. For each $N$, use the Brownian terminal reward
\begin{equation*}
h_{0,N}(s)=N^{1/3}\fh_0\left(\frac{s-N}{2N^{2/3}}\right).
\end{equation*}
Throughout this subsection, $L_N$, $\Gamma$, and $\kappa^\infty$ refer to the Brownian LPP model with this reward. Define
\begin{align}
\fh^N(x)&=N^{-1/3}L_N(2N^{2/3}x),\notag\\
\gamma_t^N(x)&=\frac{\Gamma_{\lfloor Nt\rfloor}(2N^{2/3}x)-Nt}{2N^{2/3}},
\qquad 0\le t\le1,\notag\\
\tau^N_{x,y}&=\frac{\kappa^\infty_{2N^{2/3}x,\,2N^{2/3}y}}{N}.
\label{eq:brownian rescaling}
\end{align}
As before, a pair of paths that never intersects has coalescence time $\infty$. The dependence of the Brownian environment and optimizers on $N$ is implicit; Appendix~\ref{sec:directed landscape convergence} specifies a coupling across $N$.

\begin{lemma}[Rescaled integration by parts]
\label{lem:rescaled Brownian integration by parts}
For every cylindrical functional $F$ as in \eqref{eq:cylindrical functional} and every $g_0\in C_c^1(\R)$, the quantities in \eqref{eq:brownian rescaling} satisfy \eqref{eq:rescaled Brownian integration by parts}.
\end{lemma}

\begin{proof}
Write $F(v)=G(\jbrac{g_1,v},\ldots,\jbrac{g_n,v})$. For $j=0,\ldots,n$, set
\[
f_j^N(r)=\frac{1}{2N}g_j\left(\frac{r}{2N^{2/3}}\right),
\qquad
\mathcal F_N(v)=F\!\left(x\mapsto N^{-1/3}v(2N^{2/3}x)\right).
\]
Changing variables gives
\[
\jbrac{f_j^N,L_N}=\jbrac{g_j,\fh^N},
\qquad
2N^{2/3}\jbrac{-(f_0^N)',L_N}=\jbrac{-g_0',\fh^N}.
\]
In particular, $\mathcal F_N$ is a cylindrical functional on the Brownian scale and $\mathcal F_N(L_N)=F(\fh^N)$. Its variational derivative satisfies
\begin{equation*}
D_r\mathcal F_N(L_N)=\frac{1}{2N}D_{r/(2N^{2/3})}F(\fh^N).
\end{equation*}
Apply Lemma~\ref{lem:Brownian LPP integration by parts} with $\mathcal F_N$ and $f_0^N$, and multiply the identity by $2N^{2/3}$. For $r_i=2N^{2/3}x_i$, the change of variables gives
\begin{align*}
f_0^N(r_0)\,dr_0&=N^{-1/3}g_0(x_0)\,dx_0,\\
D_{r_1}\mathcal F_N(L_N)\,dr_1&=N^{-1/3}D_{x_1}F(\fh^N)\,dx_1,\\
h_{0,N}'(\Gamma_N(r_0))&=\frac{1}{2N^{1/3}}\fh_0'(\gamma_1^N(x_0)).
\end{align*}
Also, $\sgn(r_1-r_0)=\sgn(x_1-x_0)$ and
\[
\ind{\kappa^\infty_{r_0,r_1}\le N}=\ind{\tau^N_{x_0,x_1}\le1}.
\]
Comparing the coefficients proves \eqref{eq:rescaled Brownian integration by parts}.
\end{proof}

\subsection{Directed landscape limit}

We use the Brownian terminal rewards $h_{0,N}$ and the rescaled objects $\fh^N$, $\gamma^N$, and $\tau^N$ constructed in Section~\ref{sec:rescaled Brownian LPP}. Their limiting counterparts $\fh$, $\gamma$, and $\tau$ are defined in the introduction. The normalized point-to-point value $L(r,m;s,n)$ is defined in \eqref{eq:Brownian point to point value}.

For $0\le u<v$, set
\[
m=\lfloor Nu\rfloor+1,\qquad n=\lfloor Nv\rfloor,\qquad r=Nu+2N^{2/3}x,\qquad s=Nv+2N^{2/3}y.
\]
Whenever $m\le n$ and $r\le s$, define
\[
\cL^N(x,u;y,v)=N^{-1/3}L(r,m;s,n).
\]
Applying \cite[Theorem~1.7]{dauvergne2021scaling} to Brownian last-passage percolation with $\rho=1$ and $\sigma=N^{1/3}$, and relabeling layer $k$ as $-k$, we may couple the Brownian environments with a directed landscape $\cL$ so that
\begin{equation}\label{eq:Brownian LPP field convergence}
\sup_{(x,u;y,v)\in K}\left|\cL^N(x,u;y,v)-\cL(x,u;y,v)\right|\longrightarrow 0\qquad\text{almost surely}
\end{equation}
for every compact $K\subset\{(x,u;y,v):0\le u<v\}$. All convergence statements below refer to this coupling.

\begin{lemma}[Convergence of heights, geodesics, and coalescence times]
\label{lem:landscape height path and coalescence convergence}
Under the coupling of Brownian LPP with the directed landscape described above, we have
\[
\fh^N\longrightarrow\fh\qquad\text{locally uniformly in probability}.
\]
For any two non-exceptional starting points $x$ and $y$,
\[
\sup_{t\in[0,1]}
|\gamma_t^N(x)-\gamma_t(x)|\longrightarrow 0, \qquad \tau^N_{x,y}\longrightarrow\tau_{x,y}
\]
in probability as $N\to\infty$, with coalescence times viewed in $[0,1]\cup\{\infty\}$, where $\infty$ is isolated.
\end{lemma}

\begin{proof}
\textbf{Convergence of heights and paths.} By \eqref{eq:Brownian point to line value}, the rescaled height satisfies the variational identity
\[
\fh^N(x)=\sup_{y\ge x-\frac12N^{1/3}}\left\{\cL^N(x,0;y,1)+\fh_0(y)\right\}.
\]
After reversing space and time and adjusting the normalization, \cite[Lemmas~4.4 and~4.10]{hammond2019patchwork} imply that the maximizing terminal points are tight, uniformly over starting points in compact sets. More precisely,
\[
\lim_{R \to \infty} \limsup_{N\to \infty} \P \pare{ \sup_{x\in K} |\gamma^N_1(x)| > R } = 0, \quad \text{for every compact set } K.
\]
Together with the joint field convergence \eqref{eq:Brownian LPP field convergence}, this localization implies convergence of the height.

Endpoint localization, locally uniform convergence of $\cL^N$, and uniqueness of the limiting maximizer give
\[
\gamma_1^N(x)\longrightarrow \gamma_1(x)
\]
in probability. From any subsequence, we may extract a further subsequence along which the endpoints converge almost surely. Applying \cite[Theorem~1.8]{dauvergne2021scaling} to the rescaled point-to-point optimizers with these converging endpoints gives convergence to the graph of $\gamma(x)$ in Hausdorff distance. Continuity of the limiting path turns this into uniform convergence.

\textbf{Convergence of coalescence times.} If $x=y$, then $\tau^N_{x,x}=1/N\to0=\tau_{x,x}$. By symmetry, we may therefore assume $x<y$. Write $\tau_N=\tau^N_{x,y}$ and $\tau=\tau_{x,y}$.

Set
\[
e_N:=\max_{z\in\{x,y\}}\|\gamma^N(z)-\gamma(z)\|_\infty, \qquad \omega_a = \max_{z\in\brac{x,y}}\sup_{0\le r<s\le1} \frac{|\gamma_r(z) - \gamma_s(z)|}{|r-s|^a}.
\]
\cite[Lemma~3.11, arXiv version]{ganguly2022fractal} proves that, for $a \in (0, 2/3)$
\[
\lim_{R \to \infty} \P \pare{
\omega_a > R
} = 0.
\]
We fix $a=19/30$ throughout the proof.

We first rule out premature intersections. Starting from an arbitrary subsequence, pass to a further subsequence along which $e_N\to0$ almost surely, and work on the event where this convergence holds and $\omega_a<\infty$. Let $(N_j)$ be any further subsequence, such that $\tau_{N_j}\le1$ and $\tau_{N_j}\to\tau_*\in[0,1]$. Since the optimizing paths coincide after their first intersection, we have $\gamma^{N_j}_{\tau_{N_j}}(x)=\gamma^{N_j}_{\tau_{N_j}}(y)$. Consequently,
\[
|\gamma_{\tau_*}(x)-\gamma_{\tau_*}(y)|\le2\omega_a|\tau_{N_j}-\tau_*|^a+2e_{N_j}\longrightarrow0.
\]
Thus $\tau\le\tau_*$. Every finite accumulation point along the selected subsequence is therefore at least $\tau$. If $\tau=\infty$, compactness of $[0,1]$ also implies that only finitely many of these coalescence times can be finite. Since the initial subsequence was arbitrary, the subsequence criterion gives, for every $\alpha>0$,
\[
\bP(\tau_N<\tau-\alpha,\ \tau<\infty)\longrightarrow0,\qquad \bP(\tau_N\le1,\ \tau=\infty)\longrightarrow0.
\]

Next, we rule out delayed coalescence. Fix $\alpha,\eta>0$, and set
\[
A_N=\{\tau_N>\tau+\alpha,\ \tau\le1-\eta\},\qquad q=\frac18\min\{\alpha,\eta\},\qquad \ell_N=\lfloor N q\rfloor.
\]
For integers $m\ge0$ with $m+\ell_N\le N$, define the interval
\[
I_{N,m}(z,w)=\left[m+2N^{2/3}(z-w),\,m+2N^{2/3}(z+w)\right].
\]
Let $E_{N,m}(z,z';w)$ be the event that there exist two unscaled Brownian LPP point-to-point optimizing paths that are disjoint on rows $m+1,\ldots,m+\ell_N$, with their starting points in $I_{N,m}(z,w)$ and terminal points in $I_{N,m+\ell_N}(z',w)$.

For fixed $q>0$ and $M<\infty$, \cite[Theorem~1.1]{hammond2020exponents}, applied with two polymers, implies that for all sufficiently small $w>0$ and then all sufficiently large $N$,
\begin{equation*}
\sup_{\substack{m\in\Z_{\ge0},\ m+\ell_N\le N\\ |z|,|z'|\le M+1}}\bP(E_{N,m}(z,z';w))\le C_{q,M}w^{3/2 + o(1)}.
\end{equation*}
Put $m_j=j\ell_N$ and $s_j=m_j/N$, and define, for parameters $M, \theta, w$,
\[
D^N_{M,\theta,w}=\left\{\max_{z\in\{x,y\}}\|\gamma(z)\|_\infty\le M, \; \omega_a \le \theta, \; e_N\le w/4 \right\}.
\]
On $A_N$, choose the first grid time $s_j>\tau$. These four consecutive grid times satisfy
\[
\tau<s_j<s_{j+1}<s_{j+2}<s_{j+3}\le\tau+4q<\min\{\tau+\alpha,1\}.
\]
At these times, $\gamma(x)$ and $\gamma(y)$ agree. Choose $z_0,\ldots,z_3\in w\Z$ such that
\[
|\gamma_{s_{j+i}}(x) - z_i| \le w/2, \quad i=0,1,2,3.
\]
For brevity, write $\Gamma_k^z=\Gamma_k(2N^{2/3}z)$. On $D^N_{M,\theta,w}$, for each $z\in\{x,y\}$,
\[
\left|\frac{\Gamma^z_{m_{j+i}}-m_{j+i}}{2N^{2/3}}-z_i\right|
=|\gamma^N_{s_{j+i}}(z)-z_i|\le e_N+w/2\le w,
\qquad i=1,2,3.
\]
For the starting points of the three row blocks, we have
\[
\left|\frac{\Gamma^z_{m_{j+i}+1}-m_{j+i}}{2N^{2/3}}-z_i\right|
\le e_N+\omega_aN^{-a}+\frac{1}{2N^{2/3}}+w/2\le w,
\qquad i=0,1,2,
\]
provided $N$ is sufficiently large that $\ell_N\ge2$ and $\theta N^{-a}+(2N^{2/3})^{-1}\le w/4$. Since $\tau_N>s_{j+3}$, the Brownian paths are disjoint throughout these three row blocks. On the block with rows $m+1,\ldots,m+\ell_N$, restrict each path from $(\Gamma^z_{m+1},m+1)$ to $(\Gamma^z_{m+\ell_N},m+\ell_N)$. Each restriction is a point-to-point optimizer. The estimates above place its endpoints in the required windows. Consequently,
\[
A_N\cap D^N_{M,\theta,w}\subseteq\bigcup_j\ \bigcup_{\substack{z_0,\ldots,z_3\in w\Z\\ |z_i|\le M+1}}\ \bigcap_{i=0}^2 E_{N,m_{j+i}}(z_i,z_{i+1};w).
\]
For each deterministic choice of $j,z_0,\ldots,z_3$, the three events involve disjoint sets of Brownian rows and are independent. Their intersection therefore has probability at most $C_{q,M}w^{9/2+o(1)}$. There are $O(q^{-1})$ possible time indices and $O_M(w^{-4})$ spatial tuples. The union bound gives
\[
\bP(A_N\cap D^N_{M,\theta,w})\le C_{q,M}w^{9/2+o(1)-4}\le C_{q,M}w^{1/4}.
\]
Since $e_N\to0$ in probability, for each fixed $w>0$ we obtain
\[
\limsup_{N\to\infty}\bP(A_N)
\le\bP\!\left(\max_{z\in\{x,y\}}\|\gamma(z)\|_\infty>M\right)
+\bP(\omega_a>\theta)+C_{q,M}w^{1/4}.
\]
Letting $w\downarrow0$ and then $M, \theta\uparrow\infty$ proves that $\bP(A_N)\to0$.

\begin{lemma}
\label{lem:no-terminal-atom}
For deterministic $\fh_0\in C_c^1(\R)$ and distinct deterministic
starting points $x<y$, the rightmost point-to-line geodesics satisfy
\[
\bP(\tau_{x,y}=1)=0.
\]
\end{lemma}

\begin{proof}
We first observe that, for a parabolic Airy line ensemble
$(\mathcal A_i)_{i\ge1}$ and a deterministic $g\in C^1(\R)$,
almost surely no two distinct curves $\mathcal A_i+g$ and
$\mathcal A_j+g$ have a common left-sided local maximum.
Indeed, Brownian Gibbs resampling \cite{corwin2014brownian} gives local absolute continuity
of any fixed finite collection of Airy curves with respect to
independent Brownian bridges; a deterministic $C^1$ shift preserves
this property. For independent Brownian motions $B_1,B_2$, and
a rational $r$, set
\[
R_i^r(t)=\max_{r\le s\le t}B_i(s)-B_i(t),\qquad t\ge r.
\]
By L\'evy's reflection identity, the two processes are independent
reflected Brownian motions. Simultaneous zeros at a time $t>r$
would imply a return to the origin by planar Brownian motion,
an event of probability zero. A common left-sided local maximum
would give simultaneous zeros for some rational $r$ strictly
to its left. Countable unions over rational intervals and pairs
of indices establish the observation.

Suppose $\tau_{x,y}=1$, and put
$Z=\gamma_1(x)=\gamma_1(y)$. The two geodesics are disjoint
at all interior times. By the down-the-middle isometry and its
disjoint-optimizer correspondence
\cite[Theorem~4.1(iii) and Lemma~5.2]{dauvergne2026downthemiddle},
there are two disjoint Airy switchback optimizers from the
corresponding starts to $Z$. Immediately to the left of $Z$,
their line indices are constant and distinct, say $i\ne j$.
The correspondence may be applied after a deterministic
spatial translation; localization of $Z$ and a countable
union remove the positive-quadrant restriction.

For $Z'<Z$ sufficiently close to $Z$, truncate the first
switchback optimizer at $Z'$ and jump to the top line there.
Its weight differs from the original weight by
$\mathcal A_i(Z')-\mathcal A_i(Z)$.
Since $Z$ maximizes the point-to-line objective,
\[
\mathcal A_i(Z)+\fh_0(Z)
\ge \mathcal A_i(Z')+\fh_0(Z').
\]
The same variation of the second optimizer gives
\[
\mathcal A_j(Z)+\fh_0(Z)
\ge \mathcal A_j(Z')+\fh_0(Z').
\]
Thus $\mathcal A_i+\fh_0$ and $\mathcal A_j+\fh_0$
have a common left-sided local maximum at $Z$,
contradicting the observation above.
\end{proof}

To conclude, Lemma~\ref{lem:no-terminal-atom} gives
$\bP(\tau=1)=0$.
Letting $\eta\downarrow0$ then yields
\[
\bP(\tau_N>\tau+\alpha,\ \tau<\infty)\longrightarrow0.
\]
Together with the lower bound and the assertion on $\{\tau=\infty\}$, this proves convergence of the coalescence times.
\end{proof}

\section{Tail estimates}

\begin{lemma}[Tail estimates for the KPZ fixed point]
\label{lem:landscape height tails}
Let $\fh_0\in C_c^1(\R)$. There exist $c,C>0$, depending only on $\|\fh_0\|_\infty$, such that, for $s\ge1$,
\[
\sup_{x\in\R}\bP(\fh(x)>s)\le Ce^{-cs^{3/2}},
\qquad
\sup_{x\in\R}\bP(\fh(x)<-s)\le Ce^{-cs^3}.
\]
For every $1\le q<\infty$,
\[
\sup_{x\in\R}\bE|\fh(x)|^q<\infty.
\]
\end{lemma}

\begin{proof}

Let $\fh^{\mathrm{flat}}$ denote the height function with zero reward, coupled with $\fh$ through the same directed landscape. The variational formula gives
\[
\fh^{\mathrm{flat}}(x)-\|\fh_0\|_\infty
\le \fh(x)
\le \fh^{\mathrm{flat}}(x)+\|\fh_0\|_\infty.
\]
The tail bounds follow from \cite[Proposition 4.7]{matetski2021kpz} for flat data and spatial stationarity, with constants depending only on $\|\fh_0\|_\infty$.

The moment bound follows by integrating the tails.

\end{proof}

\begin{lemma}[Tail estimates for Brownian LPP and the polymer]
\label{lem:Brownian LPP and polymer tails}
Fix $N\ge1$, a compact interval $K\subset\R$, and $h_0\in C_c^1(\R)$. There exist $c,C>0$, depending on $N,K,\|h_0\|_\infty$, such that, for $s\ge0$,
\[
\bP\left(\sup_{r\in K}|L(r)|>s\right)
+\sup_{\beta>0}\bP\left(
\sup_{r\in K}|h^\beta(r)|>s\right)
\le Ce^{-cs}.
\]
In particular, for every $1\le q<\infty$,
\[
\bE\left[\sup_{r\in K}|L(r)|^q\right]
+\sup_{\beta>0}\bE\left[
\sup_{r\in K}|h^\beta(r)|^q\right]<\infty.
\]
\end{lemma}

\begin{proof}
Let $A_K$ be the nonnegative random variable with exponential tails appearing in \eqref{eq:uniform polymer weight bounds}. Applying that estimate with $\epsilon=0$ and taking the supremum over paths gives $L(r)\le A_K$. The constant path $t_0=\cdots=t_N=r$ gives $L(r)\ge-A_K$. Hence
\[
\sup_{r\in K}|L(r)|\le A_K.
\]
Since the reference law of $w$ is a probability measure, $h^\beta(r)\le L(r)\le A_K$. On the other hand, Jensen's inequality and $\E_w[T_N^w(r)-r]=N$ give
\[
h^\beta(r)\ge\E_w[H(T^w(r))]\ge-A_K-\frac32N.
\]
Consequently, simultaneously for every $\beta>0$,
\[
\max\left\{\sup_{r\in K}|L(r)|,\sup_{r\in K}|h^\beta(r)|\right\}
\le A_K+\frac32N.
\]
The asserted tail and moment bounds follow from the exponential tail of $A_K$.
\end{proof}

\section{Flat \texorpdfstring{Airy$_1$}{Airy 1} boundary-layer estimate}
\label{app:flat-curvature}

This appendix proves the third-order boundary-layer estimate used in
Proposition~\ref{prop:flat-intensity}. 
Fix $t>0$, set $b_t=(2t)^{-1/3}$ and
\[
B_t(u,v)=b_t\operatorname{Ai}\bigl(b_t(u+v)\bigr),\qquad
Q_a=\ind{(a,\infty)},\quad
F(t,a)=\det(I-Q_aB_tQ_a).
\]
Denote the flat two-point CDF at separation $\delta>0$ by
$F_{t,\delta}(a,b)$, and put
\[
D_{t,s}(a,z)=F(t,a)-F_{t,s^2}(a,a+sz),\qquad s>0,\quad z\ge0.
\]
\begin{lemma}[Integrated Airy$_1$ boundary layer]\label{lem:flatlayer}
For $t$ in compact subsets of $(0,\infty)$,
\begin{equation}\label{eq:flat-integrated-needed}
\int_\R da\int_0^\infty dz\,D_{t,s}(a,z)
=\frac s2+s^3\partial_t\bE\fh(t,0)+o(s^3).
\end{equation}
\end{lemma}

\subsection{Airy trace-ideal bounds}
\begin{lemma}\label{flatproof:lem:trace}
Uniformly for $t$ in a compact interval of $(0,\infty)$, there is $C$ such that for $a\le0$,
\begin{equation}\label{flatproof:eq:trace}
 \|Q_aB_tQ_a\|_1\le C(1+|a|)^{7/4},\qquad
 \mathcal D_a:=\det(I+|Q_aB_tQ_a|)\le \exp\{C(1+|a|)^{7/4}\}.
\end{equation}
For $a>0$, the trace norm is $O(e^{-ca^{3/2}})$.
\end{lemma}
\begin{proof}
Set $A=1+|a|$ for $a\le0$ and partition $(a,\infty)$ into $I=[a,3A]$ and $J=[3A,\infty)$.  On $I\times I$, let $L_N=-\partial_u^2$ be the Neumann Laplacian.  For the kernel operator $T_I$,
\[
\|T_I\|_1\le \|(I+L_N)^{-1/2}\|_2\|(I+L_N)^{1/2}T_I\|_2
\le C A^{1/2}\left(\|B_t\|_{L^2(I^2)}^2+\|\partial_uB_t\|_{L^2(I^2)}^2\right)^{1/2}.
\]
The Neumann form domain is $H^1(I)$, so no endpoint vanishing is required.  Standard Airy estimates give, for $j\ge0$,
\begin{equation}\label{flatproof:eq:airy}
 |\operatorname{Ai}^{(j)}(y)|\le C_j(1+|y|)^{j/2-1/4}\quad(y\le-1),\qquad
 |\operatorname{Ai}^{(j)}(y)|\le C_j(1+y)^{j/2}e^{-cy^{3/2}}\quad(y\ge0).
\end{equation}
By changing to the sum $u+v$, one obtains
$\|B_t\|_{L^2(I^2)}^2\le CA^{3/2}$ and
$\|\partial_uB_t\|_{L^2(I^2)}^2\le CA^{5/2}$; hence $\|T_I\|_1\le CA^{7/4}$.
On the complementary rectangles $I\times J$, $J\times I$, and $J\times J$, $u+v\ge2A$.  Partition each unbounded coordinate into unit intervals and apply the same one-sided Sobolev factorization. The Airy decay in \eqref{flatproof:eq:airy} makes the sum of their trace norms bounded uniformly in $A$.  For $a>0$, the identical argument gives $O(e^{-ca^{3/2}})$.  Finally $\log\mathcal D_a\le\|Q_aB_tQ_a\|_1$.
\end{proof}

\subsection{The cofactor series and its uniform bound}
Define the smoothly extended resolvent $M_a=B_t(I-Q_aB_t)^{-1}$ and the unnormalized Fredholm minor
\[
\mathcal C_a^{(n)}(\mathbf u,\mathbf v)
=F(t,a)\det[M_a(u_i,v_j)]_{i,j=1}^n.
\]
The bordered-determinant formula, with $H_a=L^2((a,\infty))$, is
\begin{equation}\label{flatproof:eq:border}
\mathcal C_a^{(n)}(\mathbf u,\mathbf v)=
\det_{\mathbb C^n\oplus H_a}
\begin{pmatrix}
[B_t(u_i,v_j)]&-[B_t(u_i,w)]_{w>a}\\
[B_t(w,v_j)]_{w>a}&I-Q_aB_tQ_a
\end{pmatrix}.
\end{equation}
To interpret this determinant, subtract the identity on $\mathbb C^n\oplus H_a$; the resulting operator is trace class.  Formula \eqref{flatproof:eq:border} holds first for finite-rank approximants and then by trace-norm convergence.

The exact backward-heat identity for the Airy function is
\begin{equation}\label{flatproof:eq:heat}
e^{r\partial_y^2}\operatorname{Ai}(y)=e^{ry+2r^3/3}\operatorname{Ai}(y+r^2).
\end{equation}
Consequently, writing $B_{t,s}=e^{-s^2\partial_u^2}B_t$,
\begin{equation}\label{flatproof:eq:heatrow}
B_{t,s}(u,v)=b_t\exp\bigl\{-s^2b_t^3(u+v)-\tfrac23s^6b_t^6\bigr\}
\operatorname{Ai}\bigl(b_t(u+v)+s^4b_t^4\bigr).
\end{equation}
This identity is an analytic definition of the backward-heat operation on the Airy row; it is not an assertion that the backward heat semigroup is bounded on arbitrary $L^2$ functions.

Let $\mathcal C_{a,s}^{(n)}$ denote \eqref{flatproof:eq:border} with every \emph{external top row} replaced by the corresponding $B_{t,s}$ row (the interior block and external columns remain undeformed). Multilinearity gives exactly
\begin{equation}\label{flatproof:eq:row}
\mathcal C_{a,s}^{(n)}(\mathbf u,\mathbf v)
=\bigl(e^{-s^2\sum_i\partial_{u_i}^2}\mathcal C_a^{(n)}\bigr)(\mathbf u,\mathbf v).
\end{equation}

\begin{lemma}[Uniform cofactor estimate]\label{flatproof:lem:minor}
There are fixed $N,C,c>0$, uniform for $t$ in compact positive intervals, such that if $a\le0$, $A=1+|a|$, $0<s\le c A^{-N}$, $u_i=a+s\eta_i$ with $\eta_i\ge0$, and $v_i=a+s\xi_i$ with $\xi_i\le0$, then
\begin{equation}\label{flatproof:eq:cofactor}
|\mathcal C_{a,s}^{(n)}(\mathbf u,\mathbf v)|
\le\mathcal D_a\,(C\sqrt n A^N)^n
\prod_{i=1}^n (1+s|\xi_i|)^N e^{Cs^3|\xi_i|}.
\end{equation}
For each fixed total number of external derivatives, the same estimate holds with possibly larger $N$, and with factors polynomial in the relevant $|\eta_i|,|\xi_i|$. For positive $a$, after Gaussian integration in the variables $\xi_i,\eta_i$ appearing below, each fixed-order cofactor and its derivatives acquire the additional factor $e^{-ca^{3/2}}$.
\end{lemma}
\begin{proof}
The Airy bounds \eqref{flatproof:eq:airy} and explicit row formula \eqref{flatproof:eq:heatrow} imply, for each fixed $j$,
\[
\sup_{u\ge a}\|\partial_u^jB_{t,s}(u,\cdot)\mathbf1_{>a}\|_2\le CA^N,
\qquad
\|\partial_v^jB_t(\cdot,a+s\xi)\mathbf1_{>a}\|_2
\le CA^N(1+s|\xi|)^Ne^{Cs^3|\xi|}.
\]
The same bounds apply to $\partial_u^j\partial_v^kB_{t,s}(a+s\eta,a+s\xi)$. They follow by splitting $u+w\ge0$ from $u+w<0$ in the row integral: the Airy function decays in the first region and has polynomial growth in the second, whose length is at most $O(A)$; the possible exponential from \eqref{flatproof:eq:heatrow} is at most $e^{Cs^2A+Cs^3|\xi|}$, and $e^{Cs^2A}$ is bounded under the stated condition on $s$. Differentiating with respect to $s$ produces only finitely many extra factors $\eta,\xi$ and Airy derivatives.

Factor $(1+s|\xi_j|)^Ne^{Cs^3|\xi_j|}$ from the $j$th external column of the bordered determinant.  The remaining top rows, including their $n$ finite-dimensional entries, have norm at most $C\sqrt n A^N$, and the remaining external columns have norm at most $CA^N$.  Decompose the bordered determinant relative to $I_{\mathbb C^n}\oplus(I-Q_aB_tQ_a)$ into $n$ rank-one row perturbations and $n$ rank-one column perturbations.  Iterating
$|\det(I+T)|\le\det(I+|T|)$ and
$\det(I+|T+uv^*|)\le\det(I+|T|)(1+\|u\|\|v\|)$ gives \eqref{flatproof:eq:cofactor}; the harmless $1$ terms are absorbed since $A\ge1$.  External derivatives replace finitely many corresponding rows or columns, giving the same bound with fixed polynomial losses.
For positive $a$, split the Gaussian boundary variables into $s|\xi_i|\le a/2$ and its complement. In the first region every Airy argument in the cofactor lies in the positive-decay regime; in the second the Gaussian cost is $e^{-ca^2/s^2}$, which is at least as strong as the required $e^{-ca^{3/2}}$ for small $s$ and $a\ge1$. This yields the stated integrated positive-threshold bound.
\end{proof}

Let $p_1(y)=(4\pi)^{-1/2}e^{-y^2/4}$. The two-point path-integral formula of \cite[Proposition~4.3]{matetski2021kpz} and Andr\'eief's identity yield the exact, unnormalized expansion
\begin{equation}\label{flatproof:eq:series}
D_{t,s}(a,z)=\sum_{n\ge1}\frac{(-1)^{n+1}s^n}{(n!)^2}
\int_{\substack{\eta_i\ge z\\\xi_i\le0}}
\det[p_1(\eta_i-\xi_j)]_{i,j=1}^n
\mathcal C_{a,s}^{(n)}(a+s\boldsymbol\eta,a+s\boldsymbol\xi)
\,d\boldsymbol\eta\,d\boldsymbol\xi.
\end{equation}
Here the direction of the backward-heat operation is fixed by the path-integral formula \cite[Theorem~1]{quastel2013local}, and the use of unnormalized cofactors avoids separately bounding $(I-Q_aB_tQ_a)^{-1}$.

\begin{lemma}[Integrated remainder with height dependence]\label{flatproof:lem:remainder}
There exist fixed $N,C,c>0$ such that, for $a\le0$ and $0<s\le c(1+|a|)^{-N}$,
\begin{equation}\label{flatproof:eq:remainder}
\left|\int_0^\infty D_{t,s}(a,z)\,dz
-\frac{s}{2}F_a-\frac{s^2}{3\sqrt\pi}F_{aa}
-s^3\Bigl(-F_t-\frac1{12}F_{aaa}\Bigr)\right|
\le Cs^4(1+|a|)^N e^{C(1+|a|)^{7/4}}.
\end{equation}
For $a\ge0$ the right-hand side can be replaced by $Cs^4e^{-ca^{3/2}}$.
\end{lemma}
\begin{proof}
Because $\eta_i\ge z\ge0$ and $\xi_i\le0$,
\[
|\det[p_1(\eta_i-\xi_j)]|
\le n!(4\pi)^{-n/2}
\exp\left[-\tfrac14\sum_i(\eta_i^2+\xi_i^2)\right].
\]
After integrating in $z$, the domain condition contributes $\min_i\eta_i\le\sum_i\eta_i$.  Apply Lemma~\ref{flatproof:lem:minor}, absorb $e^{Cs^3|\xi|}(1+s|\xi|)^N$ into the Gaussian, and use Stirling's formula.  The integrated $n$th term is bounded by
\begin{equation}\label{flatproof:eq:n}
C\mathcal D_a\,n\,\frac{(CsA^N)^n n^{n/2}}{n!},\qquad A=1+|a|.
\end{equation}
The full sum converges uniformly when $sA^N$ is small, and its terms with $n\ge4$ sum to $O(s^4\mathcal D_a A^{N'})$.

For $n=2$, the cofactor in \eqref{flatproof:eq:series} is alternating separately in the two row and two column variables.  Twice applying the mean-value theorem, followed by the fixed-order differentiated estimate in Lemma~\ref{flatproof:lem:minor}, bounds it by $s^2|\eta_1-\eta_2||\xi_1-\xi_2|$ times the same Gaussian-integrable majorant. Its exterior prefactor is $s^2$, so its contribution is $O(s^4\mathcal D_a A^{N'})$. For $n=3$, one row subtraction already supplies an additional factor $s$, hence $O(s^4\mathcal D_a A^{N'})$.

For $n=1$, use the explicit formula \eqref{flatproof:eq:heatrow} and Taylor's formula with integral remainder at $s=0$:
\[
\begin{aligned}
B_{t,s}(a+s\eta,a+s\xi)
={}&B_t(a,a)+s(\eta\partial_1+\xi\partial_2)B_t(a,a)\\
&+s^2\bigl[(\eta^2/2-1)\partial_1^2+\eta\xi\partial_1\partial_2
 +(\xi^2/2)\partial_2^2\bigr]B_t(a,a)\\
&+O\bigl(s^3 A^{N'}(1+\eta+|\xi|)^{N'}e^{Cs^3|\xi|}\bigr).
\end{aligned}
\]
The same expansion applies to the evolved first unnormalized cofactor, since its bordered determinant depends linearly on the external Airy top row.  The Gaussian integral of the remainder is $O(s^3\mathcal D_a A^{N'})$; the exterior prefactor supplies one more $s$.

The six required Gaussian moments, obtained by putting $w=-\xi\ge0$ and integrating $\eta p_1(\eta+w)$, are
\[
\left(\tfrac12,\tfrac4{3\sqrt\pi},-\tfrac2{3\sqrt\pi},\tfrac32,-\tfrac12,\tfrac12\right)
\]
for $1,\eta,\xi,\eta^2,\eta\xi,\xi^2$ respectively. Since $B_t$ and its resolvent are symmetric, the $s^2$ coefficient is $F_{aa}/(3\sqrt\pi)$. Finally the two exact Fredholm identities
\[
F_t=-\tfrac13F(M_{11}-M_{12}+M_{22}),\qquad
F_{aaa}=F(M_{11}+2M_{12}+M_{22})
\]
show that the cubic coefficient is $-F_t-\tfrac1{12}F_{aaa}$.
Combine the first three contributions with \eqref{flatproof:eq:n} and Lemma~\ref{flatproof:lem:trace} to obtain \eqref{flatproof:eq:remainder}. The positive-$a$ bound follows by the last assertion of Lemma~\ref{flatproof:lem:minor}, including the Gaussian-tail split there.
\end{proof}

\subsection{Removal of the height cutoff}
\begin{lemma}[Integrable coefficient derivatives]\label{flatproof:lem:tails}
For each fixed compact $t$-interval, the functions $F_a,F_{aa},F_{aaa},F_t$ are integrable on $\R$, and the discarded integrals of these coefficients have cubic-exponential decay on the negative side and Airy-type decay on the positive side.
\end{lemma}
\begin{proof}
As $a\to+\infty$, the Airy row bounds and bordered Fredholm expansions give the assertion directly.  At $a=-A\to-\infty$, the GOE one-point tail gives $0\le F(t,a)\le Ce^{-cA^3}$ uniformly in $t$ on compacts.  The same differentiated-bordered-determinant argument as in Lemmas~\ref{flatproof:lem:trace}--\ref{flatproof:lem:minor} bounds $\partial_a^4F$ and $\partial_t^2F$ on unit threshold intervals and compact positive time intervals by $\exp\{CA^{7/4}\}A^N$.  Applying the elementary interpolation inequality
\[
|f^{(k)}(0)|\le C\|f\|_\infty^{1-k/m}\|f^{(m)}\|_\infty^{k/m},\qquad k<m,
\]
with $m=4$ in the height variable and $m=2$ in the time variable, shows that $F_a,F_{aa},F_{aaa},F_t$ are bounded by $Ce^{-c'A^3}$ for large $A$ (possibly with smaller $c'>0$). This gives the claimed coefficient tails.
\end{proof}
\begin{proof}[Proof of Lemma~\ref{lem:flatlayer}]
Let $L=\log(1/s)$ and truncate the height integral at
\[
 -R_-(s)\le a\le R_+(s),\qquad
 R_-(s)=(AL)^{1/3},\quad R_+(s)=(AL)^{2/3},
\]
with a sufficiently large fixed constant $A$.  Lemma~\ref{flatproof:lem:remainder} is uniform throughout the negative portion of this interval, since $s(1+R_-)^N\to0$, and its integrated error is
\[
 Cs^4(1+R_-)^{N+1}\exp\{C(1+R_-)^{7/4}\}
 =s^{4-o(1)}=o(s^3).
\]
On the positive portion its error integrates to $O(s^4)$.

The exact probabilistic identity
\[
\int_0^\infty D_{t,s}(a,z)\,dz
=\frac1s\bE\big[(Y-a)_+\mathbf1_{\{X\le a\}}\big],
\quad X=\mathfrak h(t,0),\ Y=\mathfrak h(t,s^2),
\]
combined with Cauchy--Schwarz, uniform finite height moments, and the GOE tails gives, for some $c,C>0$,
\[
\int_{\R\setminus[-R_-,R_+]}da\int_0^\infty D_{t,s}(a,z)\,dz
\le\frac Cs\big(e^{-cR_-^3}+e^{-cR_+^{3/2}}\big)=o(s^3)
\]
when $A$ is large. Lemma~\ref{flatproof:lem:tails} removes the same cutoffs from the explicit coefficients.  Integrating \eqref{flatproof:eq:remainder} now yields
\[
\begin{aligned}
\int_\R da\int_0^\infty dz\,D_{t,s}(a,z)
={}&\frac s2\int F_a\,da+\frac{s^2}{3\sqrt\pi}\int F_{aa}\,da\\
&+s^3\int\bigl(-F_t-\tfrac1{12}F_{aaa}\bigr)\,da+o(s^3)
=\frac s2+s^3m'(t)+o(s^3).
\end{aligned}
\]
Here $\int F_t\,da=-m'(t)$ follows from differentiating the tail-integral representation for $\bE\mathfrak h(t,0)$.
This is exactly \eqref{eq:flat-integrated-needed}, proving the lemma.
\end{proof}

\section*{Notation}

For convenience we collect the notation used throughout the paper.

\begingroup
\small
\renewcommand{\arraystretch}{1.15}
\begin{longtable}{@{}p{0.25\textwidth}p{0.71\textwidth}@{}}
\textbf{Symbol} & \textbf{Meaning} \\
\hline
\endfirsthead
\textbf{Symbol} & \textbf{Meaning (continued)} \\
\hline
\endhead
$\bP,\bE$ & Probability and expectation over the random environment. \\
$F(v),D_xF(v)$ & A cylindrical functional and its variational derivative at $x$. \\
\multicolumn{2}{@{}l}{\textbf{Directed landscape}} \\
$\cL(x,s;y,t)$ & The directed landscape value from $(x,s)$ to $(y,t)$, where $s<t$. \\
$\fh_0$ & The terminal reward in directed landscape coordinates. \\
$\fh,\fu$ & The point-to-line height function and its distributional derivative $\fu=\partial_x\fh$. \\
$\gamma_t(x)$ & The location at time $t$ of the rightmost optimizing point-to-line geodesic from $(x,0)$. \\
$\tau_{x,y}$ & The first meeting time of $\gamma(x)$ and $\gamma(y)$; it equals $\infty$ if they never meet. \\
$C(x),p(x)$ & For flat reward, the covariance $\cov(\fh(0),\fh(x))$ and meeting probability $\bP(\tau_{0,x}\le1)$. \\
\multicolumn{2}{@{}l}{\textbf{Brownian LPP and the polymer}} \\
$B^{(k)}$& Brownian motion on layer $k$.\\
$h_0,\beta$ & The terminal reward on the Brownian scale and the inverse temperature. \\
$\Delta_N(r,s)$ & The set of paths $r=t_0\le\cdots\le t_N=s$; $\Delta_N(r)$ allows a free terminal point. \\
$T^w(r),\E_w$ & The reference path with independent rate-one exponential increments and expectation over those increments. \\
$H(\bt)$ & The centered Brownian path weight, including the terminal reward. \\
$Z^\beta(r),h^\beta(r)$ & The point-to-line partition function and free energy. \\
$\Ppoly_r,\Epoly_r$ & Quenched polymer probability and expectation; multiple starting-point subscripts indicate conditional independence in the same environment. \\
$L(r,m;s,n)$ & The normalized point-to-point Brownian last-passage value across rows $m,\ldots,n$. \\
$L_N(r)$, or $L(r)$ & The point-to-line Brownian last-passage value with terminal reward $h_0$. \\
$\Gamma_k(r)$ & The right endpoint of the optimizing path's interval on row $k$; $\Gamma_0(r)=r$. \\
$\kappa(\bs,\bt)$ & The first row on which the two path intervals intersect, with value $\infty$ if they do not intersect. \\
$\kappa^\beta_{r_0,r_1},\kappa^\infty_{r_0,r_1}$ & The first intersection levels of two sampled polymer paths and two optimizing Brownian LPP paths, respectively. \\
\multicolumn{2}{@{}l}{\textbf{Scaling and regularity}} \\
$h_{0,N}$ & The Brownian terminal reward $N^{1/3}\fh_0((s-N)/(2N^{2/3}))$. \\
$\cL^N$ & The rescaled point-to-point Brownian last-passage field. \\
$\fh^N,\gamma^N,\tau^N$ & The rescaled Brownian height, optimizing paths, and coalescence times. \\
\hline
\end{longtable}
\endgroup

\section*{Use of AI}
We developed the main argument of this paper before generative AI tools came into widespread use. We later used these tools to help develop some of the supporting arguments and to revise the manuscript.

\section*{Related work}
The preprint \href{https://arxiv.org/abs/2611.30085}{arXiv:2611.30085} by Ruixuan Zhang overlaps with some of our results. The two papers were written independently.

\bibliographystyle{abbrvnat}
\bibliography{ref_twopoint}

@phdthesis{arratia1979coalescing,
  author       = {Arratia, Richard Alejandro},
  title        = {Coalescing {Brownian} motions on the line},
  school       = {University of Wisconsin--Madison},
  year         = {1979},
  url          = {https://search.library.wisc.edu/catalog/999512121502121}
}

@article{bakhtin2018global,
  author       = {Bakhtin, Yuri and Khanin, Konstantin},
  title        = {On global solutions of the random {Hamilton--Jacobi} equations and the {KPZ} problem},
  journal      = {Nonlinearity},
  year         = {2018},
  volume       = {31},
  number       = {4},
  pages        = {R93--R121},
  doi          = {10.1088/1361-6544/aa99a6}
}

@article{basu2023exponent,
  author       = {Basu, Riddhipratim and Busani, Ofer and Ferrari, Patrik L.},
  title        = {On the exponent governing the correlation decay of the {Airy}$_1$ process},
  journal      = {Communications in Mathematical Physics},
  year         = {2023},
  volume       = {398},
  number       = {3},
  pages        = {1171--1211},
  doi          = {10.1007/s00220-022-04544-1}
}

@article{cator2005hammersley,
  author       = {Cator, Eric and Groeneboom, Piet},
  title        = {{Hammersley}'s process with sources and sinks},
  journal      = {The Annals of Probability},
  year         = {2005},
  volume       = {33},
  number       = {3},
  pages        = {879--903},
  doi          = {10.1214/009117905000000053}
}

@misc{cator2012local,
  author       = {Cator, E. A.},
  title        = {Local maxima of two dependent {Brownian} motions never coincide},
  year         = {2012},
  howpublished = {arXiv:1207.2653},
  url          = {https://arxiv.org/abs/1207.2653}
}

@misc{dauvergne2021scaling,
  author       = {Dauvergne, Duncan and Vir{\'a}g, B{\'a}lint},
  title        = {The scaling limit of the longest increasing subsequence},
  year         = {2021},
  howpublished = {arXiv:2104.08210},
  url          = {https://arxiv.org/abs/2104.08210}
}

@article{dauvergne2022directed,
  author       = {Dauvergne, Duncan and Ortmann, Janosch and Vir{\'a}g, B{\'a}lint},
  title        = {The directed landscape},
  journal      = {Acta Mathematica},
  year         = {2022},
  volume       = {229},
  number       = {2},
  pages        = {201--285},
  doi          = {10.4310/ACTA.2022.v229.n2.a1}
}

@article{dauvergne2023last,
  author       = {Dauvergne, Duncan},
  title        = {Last passage isometries for the directed landscape},
  journal      = {Probability Theory and Related Fields},
  year         = {2023},
  volume       = {186},
  number       = {1-2},
  pages        = {391--437},
  doi          = {10.1007/s00440-022-01173-6}
}

@article{dauvergne2025geodesic,
  author       = {Dauvergne, Duncan},
  title        = {The 27 geodesic networks in the directed landscape},
  journal      = {Inventiones Mathematicae},
  year         = {2025},
  volume       = {242},
  number       = {1},
  pages        = {123--220},
  doi          = {10.1007/s00222-025-01355-8}
}

@misc{dauvergne2026downthemiddle,
  author       = {Dauvergne, Duncan and Syed, Fardin},
  title        = {Down-the-middle isometries for the {Airy} and {KPZ} sheets},
  year         = {2026},
  howpublished = {arXiv:2609.26902},
  url          = {https://arxiv.org/abs/2609.26902}
}

@article{georgiou2017geodesics,
  author       = {Georgiou, Nicos and Rassoul-Agha, Firas and Sepp{\"a}l{\"a}inen, Timo},
  title        = {Geodesics and the competition interface for the corner growth model},
  journal      = {Probability Theory and Related Fields},
  year         = {2017},
  volume       = {169},
  number       = {1-2},
  pages        = {223--255},
  doi          = {10.1007/s00440-016-0734-0}
}

@article{gu2024integration,
  author       = {Gu, Yu and Quastel, Jeremy},
  title        = {Integration by parts and invariant measure for {KPZ}},
  journal      = {Transactions of the American Mathematical Society},
  year         = {2026},
  volume       = {379},
  number       = {2},
  pages        = {1435--1452},
  doi          = {10.1090/tran/9542},
  note         = {Preprint: arXiv:2409.08465 (2024)}
}

@article{hammond2019patchwork,
  author       = {Hammond, Alan},
  title        = {A patchwork quilt sewn from {Brownian} fabric: regularity of polymer weight profiles in {Brownian} last passage percolation},
  journal      = {Forum of Mathematics, Pi},
  year         = {2019},
  volume       = {7},
  pages        = {e2},
  doi          = {10.1017/fmp.2019.2}
}

@article{hammond2020exponents,
  author       = {Hammond, Alan},
  title        = {Exponents governing the rarity of disjoint polymers in {Brownian} last passage percolation},
  journal      = {Proceedings of the London Mathematical Society},
  year         = {2020},
  volume       = {120},
  number       = {3},
  pages        = {370--433},
  doi          = {10.1112/plms.12292}
}

@article{matetski2021kpz,
  author       = {Matetski, Konstantin and Quastel, Jeremy and Remenik, Daniel},
  title        = {The {KPZ} fixed point},
  journal      = {Acta Mathematica},
  year         = {2021},
  volume       = {227},
  number       = {1},
  pages        = {115--203},
  doi          = {10.4310/ACTA.2021.v227.n1.a3}
}

@article{pimentel2016duality,
  author       = {Pimentel, Leandro P. R.},
  title        = {Duality between coalescence times and exit points in last-passage percolation models},
  journal      = {The Annals of Probability},
  year         = {2016},
  volume       = {44},
  number       = {5},
  pages        = {3187--3206},
  doi          = {10.1214/15-AOP1044}
}

@article{pimentel2022integration,
  author       = {Pimentel, Leandro P. R.},
  title        = {Integration by parts and the {KPZ} two-point function},
  journal      = {The Annals of Probability},
  year         = {2022},
  volume       = {50},
  number       = {5},
  pages        = {1755--1780},
  doi          = {10.1214/22-AOP1564}
}

@article{rahman2025infinite,
  author       = {Rahman, Mustazee and Vir{\'a}g, B{\'a}lint},
  title        = {Infinite geodesics, competition interfaces and the second class particle in the scaling limit},
  journal      = {Annales de l'Institut Henri Poincar{\'e}, Probabilit{\'e}s et Statistiques},
  year         = {2025},
  volume       = {61},
  number       = {2},
  pages        = {1075--1126},
  doi          = {10.1214/23-AIHP1454}
}

@article{rassoulagha2024coalescence,
  author       = {Rassoul-Agha, Firas and Sepp{\"a}l{\"a}inen, Timo and Shen, Xiao},
  title        = {Coalescence and total-variation distance of semi-infinite inverse-gamma polymers},
  journal      = {Journal of the London Mathematical Society},
  year         = {2024},
  volume       = {110},
  number       = {1},
  pages        = {e12955},
  doi          = {10.1112/jlms.12955}
}

@article{corwin2014brownian,
  author       = {Corwin, Ivan and Hammond, Alan},
  title        = {{Brownian} {Gibbs} property for {Airy} line ensembles},
  journal      = {Inventiones Mathematicae},
  year         = {2014},
  volume       = {195},
  number       = {2},
  pages        = {441--508},
  doi          = {10.1007/s00222-013-0462-3}
}

@misc{ganguly2022fractal,
  author       = {Ganguly, Shirshendu and Zhang, Lingfu},
  title        = {Fractal geometry of the space-time difference profile in the directed landscape via construction of geodesic local times},
  year         = {2022},
  howpublished = {arXiv:2204.01674},
  url          = {https://arxiv.org/abs/2204.01674}
}

@article{holley1975ergodic,
  author       = {Holley, Richard A. and Liggett, Thomas M.},
  title        = {Ergodic theorems for weakly interacting infinite systems and the voter model},
  journal      = {The Annals of Probability},
  year         = {1975},
  volume       = {3},
  number       = {4},
  pages        = {643--663},
  doi          = {10.1214/aop/1176996306}
}

@article{quastel2013local,
  author       = {Quastel, Jeremy and Remenik, Daniel},
  title        = {Local behavior and hitting probabilities of the {Airy}$_1$ process},
  journal      = {Probability Theory and Related Fields},
  year         = {2013},
  volume       = {157},
  number       = {3-4},
  pages        = {605--634},
  doi          = {10.1007/s00440-012-0466-8}
}

@misc{rassoulagha2026coalescence,
  author       = {Rassoul-Agha, Firas and Shen, Xiao and Zhang, Ruixuan and Zheng, Guangqu},
  title        = {Coalescence and fluctuations of {O'Connell--Yor} polymers via the free-energy correlation profile},
  year         = {2026},
  howpublished = {arXiv:2609.02064},
  url          = {https://arxiv.org/abs/2609.02064}
}
\end{document}